\documentclass[11pt]{amsart}

\usepackage{amsthm, amsfonts, amssymb, amscd}
\usepackage[pagebackref,colorlinks]{hyperref}
\usepackage{tikz-cd}
\usepackage{geometry}
\usepackage{marginnote}

\theoremstyle{definition}
\newtheorem{ntn}{Notation}[section]

\theoremstyle{plain}
\newtheorem{lem}[ntn]{Lemma}
\newtheorem{prp}[ntn]{Proposition}
\newtheorem{thm}[ntn]{Theorem}
\newtheorem{cor}[ntn]{Corollary}

\theoremstyle{definition}

\newtheorem{rem}[ntn]{Remark}
\newtheorem{exa}[ntn]{Example}

\numberwithin{equation}{section}

\newcommand{\z}{\mathbb{Z}}
\newcommand{\q}{\mathbb{Q}}

\newcommand{\F}{\mathbb{F}}

\newcommand{\EE}{\mathcal{E}}

\newcommand{\GG}{\mathcal{G}}

\newcommand{\WW}{\mathcal{W}}
\newcommand{\II}{\mathcal{I}}
\newcommand{\PP}{\mathcal{P}}
\newcommand{\RP}{\mathcal{RP}}
\newcommand{\RB}{\mathcal{RB}}
\renewcommand{\SS}{\mathcal{S}}

\newcommand{\ppp}{\mathfrak{p}}
\newcommand{\mmm}{\mathfrak{m}}

\renewcommand{\aa}{{A^\times}}
\newcommand{\tors}{{{\rm Tor}_1^{\z}}}
\newcommand{\half}{{\Big[\frac{1}{2}\Big]}}

\newcommand{\pth}{{\Big[ \frac{1}{p} \Big]}}

\newcommand{\mt}{\mapsto}
\newcommand{\lan}{\langle}
\newcommand{\ran}{\rangle}
\newcommand{\se}{\subseteq}
\newcommand{\arr}{\rightarrow}
\newcommand{\larr}{\longrightarrow}
\newcommand{\harr}{\hookrightarrow}
\newcommand{\two}{\twoheadrightarrow}
\newcommand{\Lan}{\langle \! \langle}
\newcommand{\Ran}{\rangle \! \rangle}

\newcommand{\GE}{{\rm GE}}
\newcommand{\GW}{{\rm GW}}
\newcommand{\MW}{{\rm MW}}
\newcommand{\Eb}{\mathbb{E}}
\newcommand{\Ee}{\mathrm{E}}
\newcommand{\Ind}{{\rm Ind}}

\newcommand{\Tot}{{\rm Tot}}
\newcommand{\stabe}{{\rm Stab}}
\newcommand{\GL}{\mathit{{\rm GL}}}
\newcommand{\SL}{\mathit{{\rm SL}}}

\newcommand{\SM}{\mathit{{\rm SM}}}

\newcommand{\id}{{\rm id}}
\newcommand{\im}{{\rm im}}
\newcommand{\inc}{{\rm inc}}

\renewcommand{\char}{{\rm char}}
\newcommand{\coker}{{\rm coker}}

\newcommand {\mtx}[2]
{\left(\!\!\!
	\begin{array}{cc}
		#1 \\
		#2 
	\end{array}
	\!\!\!\right)}

\newcommand{\mtxx}[4]
{\left(\!\!
	\begin{array}{cc}
		\!\!#1 & \!\!#2\\
		\!\!#3 & \!\!#4
	\end{array}\!\!
	\right)}

\begin{document}
	\title{A refined scissors congruence group and the third homology of \texorpdfstring{$\SL_2$}{Lg}}
	\author{Behrooz Mirzaii} 
	\author{Elvis Torres P\'erez}
	\address{\sf 
		Instituto de Ci\^encias Matem\'aticas e de Computa\c{c}\~ao (ICMC), Universidade de S\~ao Paulo, 
		S\~ao Carlos, Brasil}
	\email{bmirzaii@icmc.usp.br}
	\email{elvistorres@usp.br}
	
	\begin{abstract}
		There is a natural connection between the third homology of $\SL_2(A)$ and the refined Bloch 
		group $\RB(A)$ of a commutative ring $A$. In this article we investigate this connection and as 
		the main result we show that if $A$ is a universal $\GE_2$-domain such that $-1 \in \aa^2$, then 
		we have the exact sequence
		\[
		H_3(\SM_2(A),\z) \arr H_3(\SL_2(A),\z) \arr \RB(A) \arr 0,
		\]
		where $\SM_2(A)$ is the group of monomial matrices in $\SL_2(A)$. Moreover, we show that $\RP_1(A)$, 
		the refined scissors congruence  group of $A$, is naturally isomorphic to the relative homology group 
		$H_3(\SL_2(A), \SM_2(A);\z)$.
	\end{abstract}
	\maketitle
	
	For a commutative ring $A$, the study of the third homology of the group $\SL_2(A)$ is important because of its close
	connection to the third $K$-group of $A$ \cite{suslin1991}, \cite{hmm2022}, its appearance in the scissors 
	congruence problem in 3-dimensional hyperbolic geometry \cite{dupont-sah1982}, \cite{sah1989}, etc.
	
	An important method to study this group is by means of its connection to the  refined scissors congruence 
	group $\RP_1(A)$ of $A$, introduced and studied by Hutchinson \cite{hutchinson-2013}, \cite{hutchinson2013}
	\cite{hutchinson2017}, \cite{hmm2022}, \cite{C-H2022}. 
	
	Let $\RB(A)\se \RP_1(A)$ be the refined Bloch group of $A$.
	Usually there is a natural map from the third homology of $\SL_2(A)$ to  $\RB(A)$.
	In the current paper we study this map assuming minimum conditions on $A$. 
	Let $T(A)$ and $B(A)$ be the group of diagonal and upper triangular matrices in $\SL_2(A)$, 
	respectively. Assume that
	\par (i) A is a universal $\GE_2$-ring,
	\par (ii) $\mu_2(A)=\{\pm 1\}$ and $-1 \in \aa^2$,
	\par (iii) $H_n(T(A), \z) \simeq H_n(B(A), \z)$ for $n=2, 3$.\\
	Then as the first main result of this article we show that the sequence
	\begin{equation}\label{Exa1}
		H_3(\SM_2(A),\z) \arr H_3(\SL_2(A),\z) \arr \RB(A) \arr 0
	\end{equation}
	is exact, where $\SM_2(A)$ is the group of monomial matrices in $\SL_2(A)$ (Theorem~\ref{prthm}). 
	Moreover, if $A$ satisfies in conditions (i) and (iii), then as the second main result we show 
	that there is an exact sequence
	\begin{equation}\label{Exa2}
		I(A)\otimes \mu_2(A) \arr H_3(\SL_2(A), \SM_2(A); \z) \arr \displaystyle\frac{\RP_1(A)}{\lan 
			\psi_1(a^2):a\in \aa\ran} \arr 0,
	\end{equation}
	where $I(A)$ is the fundamental ideal of  $A$ (Theorem~\ref{h3-RP}). As a particular case, 
	we show that if $-1\in \aa^2$, then we have the isomorphism
	\begin{equation*}
		H_3(\SL_2(A), \SM_2(A);\z)\simeq \RP_1(A).
	\end{equation*}
	
	The homology groups of $\SL_2(A)$ relative to its subgroups $T(A)$ and $\SM_2(A)$ seems to be important.
	In this article we show that for any ring $A$ satisfying conditions (i) and (iii), we have the 
	isomorphisms
	\[
	H_2(\SL_2(A), \SM_2(A);\z)\simeq W(A), \ \ \ \  \ \ \  H_2(\SL_2(A), T(A);\z)\simeq K_1^\MW(A),
	\]
	where $W(A)$ is the Witt ring of $A$ and $K_1^\MW(A)$ is the first Milnor-Witt $K$-group of $A$.
	Moreover we show that 
	\[
	\begin{array}{c}
		H_3(\SL_2(A), T(A);\z\half) \simeq \RP_1(A)\half.
	\end{array}
	\]
	For the last two isomorphism we need to assume that $\SL_2(A)$ is perfect (see Proposition \ref{h2-h3}
	and Theorem \ref{witt}).
	
	It seems that $K_1^\MW(A)\half$ and $\RP_1(A)\half$ should be part of a chain of groups 
	\cite[App.~A]{wendt2018} with certain properties similar to  $K$-groups.
	These two groups appear in the unstable analogues of the fundamental theorem of $K$-theory
	for the second and third homology of $\SL_2$ over an infinite field \cite{hutchinson2015}, 
	which can be  used to calculate the low-dimensional homology of $\SL_2$ of Laurent polynomials 
	over certain fields. Moreover, they have certain interesting localization property 
	\cite[Theorem~6.3]{gsz2016}, \cite[Theorem~A]{hmm2022}. 
	
	We briefly outline the organization of the present paper. In Section 1 we study the $\GE_2$-rings
	and rings universal for $\GE_2$ and introduce a spectral sequence which will be our main tool 
	in handling the homology of $\SL_2(A)$. In Section 2 we introduce Hutchinson's refined scissors 
	congruence group and some of its elementary properties. In Section 3 we compare the homologies of
	the groups $T(A)$ and $B(A)$ on certain class of rings. In Section 4 we introduce the refined Bloch 
	group of $A$ and study its appearance in the spectral sequence discussed in Section 1.
	In Section 5, we study the low dimensional homologies of the group of the monomial matrices in $\SL_2(A)$.
	Section 6 is devoted to the proof of the exact sequence (\ref{Exa1}). In section 7 we develop a spectral 
	sequence suitable for the study of relative homology of groups which will be used in Section 8. 
	This spectral sequence might be known to experts but we did not find any suitable reference to the general form 
	discussed here. In Section 8 we prove the exact sequence (\ref{Exa2}) and study some other relative  homology
	groups.\\
	~\\
	{\bf Notations.} In this article all rings are commutative, except probably group rings, and have 
	the unit element $1$. For a ring $A$, $\aa$ denotes the group of invertible elements of $A$. Let
	$\WW_A$ be the set of all $a\in \aa$ such that $1-a\in \aa$. Thus
	\[
	\WW_A:=\{a\in A: a(1-a)\in \aa\}.
	\]
	Let $\GG_A:=\aa/(\aa)^2$. The element of $\GG_A$ represented by 
	$a \in \aa$ is denoted by $\lan a \ran$. We set $\Lan a\Ran:=\lan a\ran -1\in \z[\GG_A]$.\\
	~\\
	{\bf Acknowledgements.}
	The second author acknowledges that during work on this article he was supported by CAPES (Coordena\c{c}\~ao de
	Aperfei\c{c}oamento de Pessoal de N\'ivel Superior) PhD fellowship (grant number 88887.357813/2019-00).
	\section{The \texorpdfstring{$\GE_2$}{Lg}-rings and the complex of unimodular vectors}\label{GE2}
	
	Let $A$ be a commutative ring. Let  $\Ee_2(A)$ be the subgroup of $\GL_2(A)$ generated by the 
	elementary matrices
	$E_{12}(a):=\begin{pmatrix}
		1 & a\\
		0 & 1
	\end{pmatrix}$
	and
	$E_{21}(a):=\begin{pmatrix}
		1 & 0\\
		a & 1
	\end{pmatrix}$, $a\in A$.
	The group $\Ee_2(A)$ is generated by the matrices
	\[
	E(a):={\mtxx a 1 {-1} 0}, \ \ \ \ a \in  A.
	\] 
	In fact we have the following formulas
	\[
	E_{12}(a)=E(-a)E(0)^{-1}, \ \ \ \  E_{21}(a)=E(0)^{-1} E(a), \ \ \ \ 
	E(0)=E_{12}(1)E_{21}(-1)E_{12}(1).
	\]
	Let $D_2(A)$ be the subgroup of $\GL_2(A)$ generated by diagonal matrices. 
	Let $\GE_2(A)$ be the subgroup of $\GL_2(A)$  generated 
	by $D_2(A)$ and $\Ee_2(A)$. A ring $A$ is called a {\it $\GE_2$-ring} if 
	\[
	\GE_2(A)=\GL_2(A).
	\]
	Since $\Ee_2(A)=\SL_2(A)\cap \GE_2(A)$ and $\GL_2(A)=\SL_2(A)D_2(A)$, 
	this condition is equivalent to $\Ee_2(A)=\SL_2(A)$.
	
	For any $a\in \aa$, let
	$D(a):=\begin{pmatrix}
		a & 0\\
		0 & a^{-1}
	\end{pmatrix}$. Observe that $D(-a)=E(a)E(a^{-1})E(a)$. Thus $D(a) \in \Ee_2(A)$.
	For any $x,y\in A$ and $a\in \aa$, we have the following relations between matrices 
	$E(x)$ and $D(a)$: 
	
	\begin{itemize}
		\item[(1)] $E(x)E(0)E(y)=D(-1)E(x+y)$,
		\item[(2)] $E(x)D(a)=D(a^{-1})E(a^2x)$,
		\item[(3)] $D(a)D(b)=D(ab)$.
	\end{itemize}
	
	A ring $A$ is {\it  universal for} $\GE_2$ if $\Ee_2(A)$ is generated by the elements $E(x)$ 
	subject to relations (1)-(3), where
	\[
	D(a):= E(-a)E(-a^{-1})E(-a).
	\]
	A $\GE_2$-ring which is universal for $\GE_2$ is called  a {\it universal $\GE_2$-ring}. 
	Thus a universal $\GE_2$-ring is characterized by the property that $\SL_2(A)$ is generated 
	by the matrices $E(x)$ with (1)-(3) as a complete set of defining relations, where $D(a)$ is defined in above.
	
	Any local ring is a universal $\GE_2$-ring \cite[Theorem 4.1]{cohn1966}.
	Moreover, Euclidean domains are $\GE_2$-rings \cite[\S 2]{cohn1966}. For more example of 
	$\GE_2$-rings and rings universal for $\GE_2$ see \cite{cohn1966} and \cite{hutchinson2022}.
	
	A (column) vector ${\pmb u}={\mtx {u_1} {u_2}}\in A^2$ is said to be unimodular 
	if there exists a vector ${\pmb v}={\mtx {v_1} {v_2}}$ such that the matrix 
	$({\pmb u}, {\pmb v}):={\mtxx {u_1} {v_1} {u_2} {v_2}}$ is  an invertible matrix. 
	
	For any non-negative integer $n$, let $X_n(A^2)$ be the free abelian group generated by 
	the set of all $(n+1)$-tuples $(\lan{\pmb v_0}\ran, \dots, \lan{\pmb v_n}\ran)$, where 
	every ${\pmb v_i} \in A^2$ is unimodular and for any two distinct vectors 
	${\pmb v_i}, {\pmb v_j}$, the matrix ${\pmb v_i}, {\pmb v_j}$ is invertible. 
	Observe that $\lan{\pmb v}\ran\se A^2$ is the line $\{{\pmb v}a: a\in A\}$.
	
	We consider $X_{l}(A^2)$ as a left $\GL_2(A)$-module (resp. left $\SL_2(A)$-module) in 
	a natural way. If necessary, we convert this action to a right action by the definition
	$m.g:=g^{-1}m$. Let us define the $l$-th differential operator
	\[
	\partial_l : X_l(A^2) \arr X_{l-1}(A^2), \ \ l\ge 1,
	\]
	as an alternating sum of face operators which throws away the $i$-th component of 
	generators. Let $\partial_{-1}=\epsilon: X_0(A^2) \arr \z$ be defined by 
	$\sum_i n_i(\lan v_{0,i}\ran) \mt \sum_i n_i$. Hence we have the complex
	\[
	X_\bullet(A^2)\arr \z: \  \cdots \larr  X_2(A^2) \overset{\partial_2}{\larr} X_1(A^2) 
	\overset{\partial_1}{\larr} X_0(A^2) \arr \z \arr 0.
	\]
	We say that the above complex is exact in dimension $<k$ if the complex
	\[
	X_k(A^2)\overset{\partial_k}{\larr} X_{k-1}(A^2)\overset{\partial_{k-1}}{\larr} \cdots 
	\overset{\partial_2}{\larr} X_1(A^2) \overset{\partial_1}{\larr} X_0(A^2) \arr \z \arr 0
	\]
	is exact.
	
	\begin{prp}[Hutchinson]\label{GE2C}
		Let $A$ be a commutative ring.
		\par {\rm (i)} The ring $A$ satisfies the condition that $X_\bullet(A^2)  \arr \z$ is 
		exact in dimension $<1$ if and only if $A$ is a $\GE_2$-ring. 
		\par {\rm (ii)} If $A$ is universal for  $\GE_2$, then  $X_\bullet(A^2)$ is exact in 
		dimension $1$, i.e. $H_1(X_\bullet(A^2))=0$.
	\end{prp} 
	\begin{proof}
		See \cite[Theorem 3.3, Theorem 7.2 and Corollary 7.3]{hutchinson2022}.
	\end{proof}
	
	\begin{rem}
		In \cite[Theorem 3.3, Theorem 7.2]{hutchinson2022}  Hutchinson calculated $H_0$ and $H_1$ 
		of the complex $X_\bullet(A^2)$ for any commutative ring $A$. 
	\end{rem}
	
	Assume that A satisfies the condition that $X_\bullet(A^2)\arr \z$ is exact in dimension $<1$, (i.e. $A$ is a $\GE_2$-ring 
	by Proposition \ref{GE2C}). Let 
	\[
	Z_1(A^2):=\ker(\partial_1).
	\]
	From the complex
	\begin{equation}\label{comp1}
		0 \arr Z_1(A^2) \overset{\inc}{\arr} X_1(A^2) \overset{\partial_1}{\arr} X_0(A^2) \arr 0,
	\end{equation}
	we obtain the double complex
	\[
	D_{\bullet,\bullet}: 0 \arr F_\bullet \otimes_{\SL_2(A)} Z_1(A^2) 
	\overset{\id_{F_\bullet}\otimes\inc}{-\!\!\!-\!\!\!-\!\!\!\larr} 
	F_\bullet \otimes_{\SL_2(A)}X_1(A^2)  \overset{\id_{F_\bullet}\otimes\partial_1}{-\!\!\!-\!\!\!-\!\!\!\larr} 
	F_\bullet \otimes_{\SL_2(A)}X_0(A^2)  \arr 0,
	\]
	where $F_\bullet \arr \z$ is a projective resolution of $\z$ over $\SL_2(A)$.
	This gives us the first quadrant spectral sequence
	\[
	E^1_{p.q}=\left\{\begin{array}{ll}
		H_q(\SL_2(A),X_p(A^2)) & p=0,1\\
		H_q(\SL_2(A),Z_1(A^2)) & p=2\\
		0 & p>2
	\end{array}
	\right.
	\Longrightarrow H_{p+q}(\SL_2(A),\z).
	\]
	In our calculations  we usually use  the bar resolution $B_\bullet(\SL_2(A))\arr \z$ 
	\cite[Chap.I, \S 5]{brown1994}. 
	
	The group $\SL_2(A)$ acts transitively on the set of generators of $X_i(A^2)$ for $i=0,1$. Let
	\[
	{\pmb\infty}:=\lan {\pmb e_1}\ran, \ \ \  {\pmb 0}:=\lan {\pmb e_2}\ran , \ \ \  
	{\pmb a}:=\lan {\pmb e_1}+ a{\pmb e_2}\ran, \ \ \ a\in \aa,
	\]
	where ${\pmb e_1}:={\mtx 1 0}$ and $ {\pmb e_2}:={\mtx 0 1}$.
	We choose $({\pmb \infty})$ and $({\pmb \infty} ,{\pmb 0})$ as 
	representatives of the orbit of the generators of $X_0(A^2)$ and $X_1(A^2)$,
	respectively.  Therefore 
	\[
	X_0(A^2)\simeq \Ind _{B(A)}^{\SL_2(A)}\z, \ \ \ \ \ \ \ \ \ \ \ \ 
	X_1(A^2)\simeq \Ind _{T(A)}^{\SL_2(A)}\z,
	\]
	where 
	\[
	B(A):=\stabe_{\SL_2(A)}({\pmb\infty})=\bigg\{\begin{pmatrix}
		a & b\\
		0 & a^{-1}
	\end{pmatrix}:a\in \aa, b\in A\bigg\},
	\]
	\[
	T(A):=\stabe_{\SL_2(A)}({\pmb \infty},{\pmb 0})=\bigg\{\begin{pmatrix}
		a & 0\\
		0 & a^{-1}
	\end{pmatrix}:a\in \aa\bigg\}.
	\]
	Note that $T(A)\simeq \aa$. In our calculations usually we identify $T(A)$ 
	with $\aa$. Thus by Shapiro's lemma we have
	\[
	E_{0,q}^1 \simeq H_q(B(A),\z), \ \ \ \ \ \ E_{1,q}^1 \simeq H_q(T(A),\z).
	\]
	In particular, $E_{0,0}^1\simeq \z\simeq E_{1,0}^1$. Moreover
	$d_{1, q}^1=H_q(\sigma) - H_q(\inc)$,
	where $\sigma: T(A) \arr B(A)$ is given by $\sigma(X)= wX w^{-1}=X^{-1}$ for 
	$w:=E(0)={\mtxx 0 1 {-1} 0}$. This easily implies that  $d_{1,0}^1$ is trivial, 
	$d_{1,1}^1$  is induced by the map $T(A) \arr B(A)$, $X\mt X^{-2}$, and 
	$d_{1,2}^1$ is trivial.  Thus $\ker(d_{1,1}^1)=\mu_2(A)=\{a\in \aa: a^2=1\}$.
	It is straightforward to check that $d_{2,0}^1:H_0(\SL_2(A),Z_1(A^2)) \arr \z$ is surjective and for any 
	$b\in \mu_2(A)$, $d_{2,1}^1([b]\otimes \partial_2({\pmb \infty}, {\pmb 0}, {\pmb a}))=b$. 
	Hence $E_{1,0}^2=0$ and $E_{1,1}^2=0$.
	
	\section{The refined scissors congruence group}
	
	Let $Z_2(A^2):=\ker(\partial_2)$. Following Coronado and Hutchinson \cite[\S3]{C-H2022} we define
	\[
	\RP(A):=H_0(\SL_2(A), Z_2(A^2))=Z_2(A^2)_{\SL_2(A)}.
	\]
	Note that $\RP(A)$ is a $\GG_A$-module.  The inclusion $\inc: Z_2(A^2) \arr X_2(A^2)$ induces the map
	\[
	\lambda: \RP(A)=Z_2(A^2)_{\SL_2(A)}  \overset{\overline{\inc}}{\larr} X_2(A^2)_{\SL_2(A)}.
	\]
	The orbits of the action of  $\SL_2(A)$ on $X_2(A)$ is represented by
	$\lan a\ran[\ ]:=({\pmb\infty}, {\pmb 0},{\pmb a})$, $\lan a\ran\in \GG_A$.
	Therefore $X_2(A^2)_{\SL_2(A)} \simeq \z[\GG_A]$. The $\GG_A$-module  
	\[
	\RP_1(A):=\ker\Big(\lambda:\RP(A) \arr \z[\GG_A]\Big)
	\]
	is called the {\it refined scissors congruence group} of $A$.  We call 
	\[
	\GW(A):=H_0(\SL_2(A), Z_1(A^2))
	\]
	the {\it Grothendieck-Witt group} of $A$ (see Remark \ref{GW} for an explanation on the choice of this name). 
	Let $\epsilon:=d_{2,0}^1: \GW(A) \arr \z$.
	The kernel of $\epsilon$ is called the {\it fundamental ideal} of $A$ and is denoted by $I(A)$.
	
	Consider the sequence
	\[
	X_4(A^2)_{\SL_2(A)} \overset{\overline{\partial_4}}{\arr} X_3(A^2)_{\SL_2(A)}
	\overset{\overline{\partial_3}}{\arr} \RP(A) \arr 0
	\]
	of $\GG_A$-modules. The orbits of the action of  $\SL_2(A)$ on $X_3(A)$  and $X_4(A)$ 
	are  represented by
	\[
	\lan a\ran[x]:=({\pmb\infty}, {\pmb 0},{\pmb a}, \pmb{ax}),
	\ \  \text{and}  \ \   
	\lan a\ran[x,y]:= ({\pmb\infty}, {\pmb 0},{\pmb a}, \pmb{ax}, \pmb{ay}),
	\ \ \lan a\ran\in \GG_A, x,y,x/y\in \WW_A,
	\]
	respectively. Thus $X_3(A^2)_{\SL_2(A)} $ is the free $\z[\GG_A]$-module generated by 
	the symbols $[x]$, $x\in \WW_A$ and  $X_4(A^2)_{\SL_2(A)} $ is the free 
	$\z[\GG_A]$-module generated by the symbols
	$[x,y]$, $x,y,x/y\in \WW_A$. It is straightforward to check that 
	\[
	\overline{\partial_4}([x,y])=[x]-[y]+\lan x\ran\bigg[\frac{y}{x}\bigg]-
	\lan x^{-1}-1\ran\Bigg[\frac{1-x^{-1}}{1-y^{-1}}\Bigg] + 
	\lan 1-x\ran\Bigg[\frac{1-x}{1-y}\Bigg].
	\]
	Let $\overline{\RP}(A)$ be the quotient of the free $\GG_A$-module generated 
	by the symbols $[x]$, $x\in \WW_A$ over the subgroup generated by the elements
	\[
	[x]-[y]+\lan x\ran\bigg[\frac{y}{x}\bigg]-
	\lan x^{-1}-1\ran\Bigg[\frac{1-x^{-1}}{1-y^{-1}}\Bigg] + 
	\lan 1-x\ran\Bigg[\frac{1-x}{1-y}\Bigg],  \ \ \ \ \ \ x,y,x/y\in \WW_A.
	\]
	Clearly we have the natural map $\overline{\RP}(A) \arr \RP(A)$. It is straightforward 
	to check that the composite
	\[
	\overline{\RP}(A) \arr \RP(A) \overset{\lambda}{\larr} \z[\GG_A]
	\]
	is given by
	\[
	[x]\mapsto -\Lan x\Ran\Lan 1-x\Ran.
	\]
	Let $\overline{\RP}_1(A)$ be the kernel of this composite. Thus we have a natural map
	\[
	\overline{\RP}_1(A) \arr \RP_1(A).
	\]
	The sequence
	\[
	X_3(A^2)_{\SL_2(A)} \overset{\overline{\partial_3}}{\arr} X_2(A^2)_{\SL_2(A)}    
	\overset{\overline{\partial_2}}{\arr} \GW(A) \arr 0
	\]
	induces the natural map 
	\[
	\overline{\GW}(A):=\z[\GG_A]/\Big\lan \Lan a\Ran \Lan 1-a\Ran: a\in \WW_A\Big\ran \arr \GW(A).
	\]
	Let $\II_A$ be the kernel of the augmentation map $\z[\GG_A] \arr \z$ and set 
	\[
	\overline{I}(A):=\II_A/\Big\lan \Lan a\Ran \Lan 1-a\Ran: a \in \WW_A\Big\ran.
	\]
	Thus we have a natural map $\overline{I}(A)\arr I(A)$.
	
	If the complex  $X_\bullet(A^2)\arr \z$ is exact in dimension $<2$, then 
	$\overline{I}(A)\arr I(A)$ is surjective. If the complex is exact in dimension 
	$<3$, then the maps 
	\[
	\text{$\overline{\RP}(A) \arr \RP(A)$ and $\overline{\RP}_1(A) \arr \RP_1(A)$}
	\]
	are surjective and $\overline{I}(A)\simeq I(A)$. Moreover, if the complex is exact 
	in dimension $<4$, then  $\overline{\RP}(A) \simeq \RP(A)$ and 
	$\overline{\RP}_1(A) \simeq \RP_1(A)$.
	
	\begin{rem}\label{II-I}
		Assume that A satisfies the condition that $X_\bullet(A^2)\arr \z$ is exact in dimension $<2$. From the exact sequence
		\[
		0 \arr  Z_2(A^2) \arr X_2(A^2) \arr Z_1(A^2) \arr 0
		\]
		we obtain the exact sequence $\RP(A) \overset{\lambda}{\larr} \z[\GG_A] \arr \GW(A) \arr 0$.
		This induces the exact sequence 
		\[
		\RP(A) \overset{\lambda}{\larr} \II_A \arr I(A) \arr 0.
		\]
		For any $a\in\aa$, let
		\[
		\psi_1(a):=({\pmb \infty}, {\pmb 0}, {\pmb a})+({\pmb 0}, {\pmb \infty}, {\pmb a})
		-({\pmb \infty}, {\pmb 0}, {\pmb 1})-({\pmb 0}, {\pmb \infty}, {\pmb 1}) \in \RP(A).
		\]
		Then 
		\[
		\lambda(\psi_1(a))=p_{-1}^+\Lan a\Ran,
		\]
		where $p_{-1}^+:=\lan-1\ran +1\in \z[\GG_A]$. This induces a natural surjection 
		\[
		\II_A/p_{-1}^+\II_A \two I(A).
		\]
		We borrowed this general form of the element $\psi_1(a)$ from \cite{C-H2022}. For more on this element see \cite[\S3.2]{C-H2022},
		\cite{hutchinson2013}, \cite{hutchinson2017} and \cite{hmm2022}.
	\end{rem}
	
	\begin{rem}\label{GW}
		Assume that $A$ satisfies the condition that $X_\bullet(A^2)\arr \z$ is exact in dimension $< 3$.
		Then, as we have seen in above, the Grothendieck-Witt group $\GW(A)$ is isomorphic to 
		\[
		\overline{\GW}(A)=\z[\GG_A]/\lan \Lan a\Ran\Lan 1-a\Ran:a\in \WW_A\ran.
		\]
		When $A$ is a field of characteristic not $2$, or a local ring with sufficiently
		large residue field $k$, where $\char(k)\neq 2$, this latter group coincides with the classical Grothendieck-Witt 
		group of symmetric bilinear forms over $A$ (see \cite[Lemma 1.1, Chap. 4]{HM1973} or 
		\cite[Theorems 4.1,4.3, Chap. II]{lam2005}).
		When $A$ is a more general local ring this result is also well-known and for a proof of this
		see \cite[Lemma 1.7, Theorem 4.1]{gsz2016}.
	\end{rem}
	
	\section{ The map \texorpdfstring{$H_n(T(A),\z)\arr H_n(B(A),\z)$}{Lg} } \label{iso}
	The groups $B(A)$ and T(A) sit in the extension $1\arr N(A) \arr B(A) \arr T(A) \arr 1$, where
	\[
	N(A):=\Bigg\{\begin{pmatrix}
		1 & b\\
		0 & 1
	\end{pmatrix}: b\in A \Bigg \}\simeq A.
	\]
	This extension splits canonically and $T(A)$ acts as follow on $N$: 
	\[
	a.\begin{pmatrix}
		1 & b\\
		0 & 1
	\end{pmatrix} :=  \begin{pmatrix}
		a & 0\\
		0 & a^{-1}
	\end{pmatrix}  \begin{pmatrix}
		1 & b\\
		0 & 1
	\end{pmatrix} 
	{\begin{pmatrix}
			a & 0\\
			0 & a^{-1}
	\end{pmatrix}}^{-1}=\begin{pmatrix}
		1 & a^2b\\
		0 & 1
	\end{pmatrix}.
	\]
	So if we assume that $T(A)=\aa$ and $N(A)=A$, then the action of $\aa$ on $A$ is 
	given by $a.b:=a^2b$. Thus
	\[
	H_n(B(A),\z) \simeq H_n(T(A),\z) \oplus H_n(B(A),T(A);\z).
	\]
	We denote the relative homology group $H_n(B(A),T(A);\z)$ by $\SS_n$.
	(See Section \ref{RP-SL} for an exact sequence involving this relative homology group).
	
	By studying the Lyndon/Hochschild-Serre spectral sequence of the above extension, 
	it follows that
	\[
	\SS_1\simeq H_0(\aa,A)=A_\aa=A/\lan a^2-1|a\in \aa \ran
	\]
	and $\SS_2$ sits in the exact sequence 
	\[
	H_2(\aa,  A)  \arr H_2(A,\z)_\aa \arr \SS_2 \arr H_1(\aa, A)\arr 0.
	\]
	
	\begin{lem}\label{lem-sus}
		Let $G$ be an abelian group, $A$ a commutative ring, $M$ an $A$-module and 
		$\varphi: G \arr A^\times$ a homomorphism of groups which turns $A$ and $M$ 
		into $G$-modules.  If $H_0(G,A)=0$, then for any  $n\geq 0$, $H_n(G,M)=0$.
	\end{lem}
	\begin{proof}
		See \cite[Lemma~1.8]{nes-suslin1990}.
	\end{proof}
	
	\begin{cor}\label{H0A}
		Let $A$ be a ring and let $\aa$ acts on $A$ as $a.x:= a^2x$. If $H_0(\aa,A)=0$, 
		then $H_n(\aa,A)=0$ for any $n\geq 0$.
	\end{cor}
	\begin{proof}
		Use the above lemma by considering $\varphi: \aa \arr \aa$, $a\mt a^2$.
	\end{proof}
	
	\begin{exa}\label{H00A}
		(i) If $A$ is a local ring such that $|A/\mmm_A|>3$, then always we can find $a\in \aa$
		such that $a^2-1\in \aa$. Thus $H_0(\aa,A)=0$.
		\par (ii) Let $A$ be a ring such that  $6\in \aa$. Then 
		\[
		1=3(2^2-1)+(-1)(3^2-1) \in \lan a^2-1:a\in \aa\ran. 
		\]
		Hence $H_0(\aa,A)=0$.
	\end{exa}
	
	\begin{exa}
		If $H_0(\aa,A)=0$, then by the above corollary $H_n(\aa,A)=0$ for $n\geq 0$. Thus $\SS_1=0$ 
		and $\SS_2\simeq H_2(A,\z)_\aa$. Therefore $H_1(T(A),\z)\simeq H_1(B(A),\z)$ and we have the 
		exact sequence
		\[
		0 \arr H_2(A,\z)_\aa \arr H_2(B(A),\z)\arr H_2(T(A),\z) \arr 0.
		\]
		Moreover, we have the exact sequence
		\[
		H_3(A,\z)_\aa \arr \SS_3 \arr H_1(\aa, A\wedge A) \arr 0.
		\]
	\end{exa}

	\begin{lem}\label{1/6}
		If $A$ is a subring of $\q$, then for any $n\geq 0$, 
		\[
		H_n(B(A),\z)\simeq H_n(T(A),\z) \oplus H_{n-1}(\aa, A).
		\]
		In particular if $6\in \aa$, then $H_n(T(A),\z) \simeq H_n(B(A),\z)$.
	\end{lem}
	\begin{proof}
		It is well known that any finitely generated subgroup of $\q$ is cyclic. Thus
		$A$ is a direct limit of infinite cyclic groups. Since $H_n(\z,\z)=0$ for any
		$n\geq  2$ \cite[page 58]{brown1994} and since homology commutes with direct limit 
		\cite[Exer.~6, \S~5, Chap.~V]{brown1994}, we have $H_n(A,\z)=0$ for $n\geq 2$. Now 
		the claim follows from an easy analysis of the Lyndon/Hochschild-Serre spectral 
		sequence associated to the split extension $1\arr  N(A) \arr B(A) \arr T(A) \arr 1$.
		
		If  $6\in \aa$,  then by Example~\ref{H00A}(ii) we have $H_0(\aa,A)=0$. So by Corollary~\ref{H0A},
		$H_n(\aa, A)=0$ for any $n$. Therefore the claim follows from the first part of the lemma.
	\end{proof}
	
	\begin{exa}\label{HnB}
		(i) Let $A=\z$. Since $\z^\times=\{\pm 1\}$, the action of $\z^\times$ on $A=\z$
		is trivial. Thus $H_n(\z^\times, \z)$ is $\z$ if $n=0$, is trivial if $n$ is
		even and is $\z/2$ if $n$ is odd. Now by the previous lemma we have
		\[
		H_1(B(\z),\z)\simeq H_1(T(\z),\z)\oplus \z,
		\]
		and for any positive integer $m$,
		\[
		H_{2m}(B(\z),\z)\simeq H_{2m}(T(\z),\z) \oplus \z/2\simeq \z/2, 
		\]
		\[
		H_{2m+1}(B(\z),\z)\simeq H_{2m+1}(T(\z),\z) \simeq \z/2.
		\]
		\par (ii) Let $p$ be a prime and let $A:=\z_{(p)}=\{a/b\in \q| a, b\in \z, p\nmid b\}$. 
		Then $\z_{(p)}$ is local and its residue field  is isomorphic to $\F_p$. If $p\neq 2, 3$, 
		then the residue field of $A$ has more than $3$ elements. Thus 
		\[
		H_n(T(\z_{(p)}),\z) \simeq H_n(B(\z_{(p)}),\z)
		\]
		for any $n\geq 0$ (Example \ref{H00A}). 
		
		Let $B=\z_{(2)}$. Consider the action of $B^\times$ on $\q$ as usual: $b.x:=b^2x$. It is 
		straightforward to check that $H_0(B^\times, \q)=0$. Thus by Lemma \ref{lem-sus},
		$H_n(B^\times, \q)=0$ for any $n\geq 0$. 
		Consider the exact sequence $0 \arr B \arr \q \arr \q/B \arr 0$. Note that 
		$\q/B\simeq \z_{2^\infty}:=\z\half/\z$. From the long exact sequence associated to this short exact 
		sequence, we obtain
		\[
		H_{n-1}(B^\times, B) \simeq H_{n}(B^\times, \z_{2^\infty}), \ \ \ \ \ \   n\geq 1.
		\]
		We have a similar result for $B=\z_{(3)}$. Therefore for $p=2,3$, we have
		\[
		H_n(B(\z_{(p)}),\z)\simeq H_n(T(\z_{(p)}),\z) \oplus H_{n}(\z_{(p)}^\times, \z_{p^\infty}).
		\]
		Note that $H_{n}(\z_{(2)}^\times, \z_{2^\infty})$ and $H_{n}(\z_{(3)}^\times, \z_{3^\infty})$ 
		are 2-power and 3-power torsion groups, respectively. One easily can show that
		$H_0(\z_{(2)}^\times , \z_{(2)})\simeq \z/8$ and $H_0(\z_{(3)}^\times , \z_{(3)})\simeq \z/3$.
	\end{exa}
	
	\begin{lem} \label{1/2}
		Let $p$ be a prime number and let $A_p=\z[\frac{1}{p}]$. Then 
		\par {\rm (i)} $H_1(B(A_p),\z)\simeq H_1(T(A_p),\z)\oplus \z/(p^2-1)$,
		\par {\rm (ii)} for any $n\geq 2$, $H_n(T(A_2),\z)\simeq H_n(B(A_2),\z)$,
		\par {\rm (iii)} for $p\neq 2$ and $n\geq 2$, we have $H_n(B(A_p),\z) 
		\simeq H_n(T(A_p),\z) \oplus \z/2$.
	\end{lem}
	\begin{proof}
		We need to calculate $H_n(A_p^\times, A_p)$. The rest follows from Lemma \ref{1/6}.
		In the following we will use the calculation of the homology groups of cyclic groups
		\cite[page 58]{brown1994}. 
		
		From the extension $1 \arr \mu_2(A_p) \arr A_p^\times \arr \lan p\ran \arr 1$
		we obtain the Lyndon/Hochschild-Serre spectral sequence
		\[
		{E'}_{r,s}^2=H_r(\lan p \ran, H_s(\mu_2(A_p), A_p))\Rightarrow H_{r+s}(A_p^\times, A_p).
		\]
		Since $\lan p\ran$  is an infinite cyclic group, we have ${E'}_{r,s}^2=0$ for $r\geq 2$. 
		Moreover
		\[
		H_s(\mu_2(A_p), A_p)\simeq
		\begin{cases}
			A_p & \text{if $s=0$}\\
			A_p/2 & \text{if $s$ is odd}\\
			0 & \text{if $s$ is even.}
		\end{cases}
		\]
		\par (i) This follows from the isomorphism
		$H_0(A_p^\times, A_p)= A_p/\lan p^2-1\ran\simeq \z/(p^2-1)$. 
		\par (ii) Since $2\in A_2^\times$, $A_2/2=0$. This implies that ${E'}_{r,s}^2=0$ for any
		$s\geq 1$. Now from the above spectral sequence we obtain $H_n(A_2^\times, A_2)=0$ for any 
		$n\geq 1$.
		\par (iii) We need to calculate ${E'}_{0,s}^2$ and ${E'}_{1,s}^2$ for any $s\geq 1$. 
		Note that $A_p/2\simeq \z/2$. Now it is easy to see that $H_0(\lan p\ran, A_p/2)\simeq\z/2$ 
		and $H_1(\lan p\ran, A_p/2)\simeq\z/2$. Thus for any $s\geq 1$,
		\[
		{E'}_{0,s}^2\simeq {E'}_{1,s}^2\simeq
		\begin{cases}
			0 & \text{if $s$ is even}\\
			\z/2 & \text{if $s$ is odd.}
		\end{cases}
		\]
		Now from the above spectral sequence it follows that $H_n(A_p^\times, A_p)\simeq \z/2$ for any $n\geq 1$.
	\end{proof}
	
	\begin{prp}\label{iso-hut}
		\par {\rm (i)} Let $A$ be a local domain such that either $A/\mmm_A$ is infinite 
		or if $|A/\mmm_A|=p^d$, we have $(p -1)d > 2n$. Then $H_n(T(A), \z)\simeq H_n(B(A),\z)$.
		\par {\rm (ii)} Let $A$ be a local ring such that either $A/\mmm_A$ is infinite 
		or if $|A/\mmm_A|=p^d$, we have $(p -1)d > 2(n+1)$. Then $H_n(T(A), \z)\simeq H_n(B(A),\z)$.
	\end{prp}
	\begin{proof}
		(i) For this see \cite[Proposition 3.19]{hutchinson2017}.
		\par (ii) Similar to the proof of part (i) presented in \cite[Proposition 3.19]{hutchinson2017},
		we can show that $H_n(T(A), k)\simeq H_n(B(A), k)$, where $k$ is a prime field and $(p -1)d > 2n$.
		Now the claim follows from \cite[Lemma 2.3]{mirzaii2017}.
	\end{proof}
	
	\section{The refined Bloch group}
	
	Assume that $A$ satisfies the condition that $X_\bullet(A^2)\arr \z$ is exact in dimension $< 2$.
	Then from the exact sequence
	\[
	0 \arr  Z_2(A^2) \arr X_2(A^2) \arr Z_1(A^2) \arr 0
	\]
	we obtain the  long exact sequence
	\[
	H_1(\SL_2(A), Z_2(A^2))\! \arr \!H_1(\SL_2(A),X_2(A^2))\! \arr\! H_1(\SL_2(A),Z_1(A^2))
	\!\overset{\delta}{\arr}\! H_0(\SL_2(A),Z_2(A^2)) 
	\]
	\[
	\arr H_0(\SL_2(A),X_2(A^2)) \arr H_0(\SL_2(A),Z_1(A^2)) \arr 0.
	\]
	Choose $({\pmb \infty}, {\pmb 0}, {\pmb a})$, $\lan a\ran \in \GG_A$, as representatives 
	of the orbits of the generators of $X_2(A^2)$. Then
	\[
	X_2(A^2)\simeq \bigoplus_{\lan a\ran\in \GG_A} \Ind _{\mu_2(A)}^{\SL_2(A)}\z\lan a\ran, 
	\]
	where $\mu_2(A)\simeq\stabe_{\SL_2(A)}({\pmb \infty},{\pmb 0}, {\pmb a})$. Thus
	\[
	H_1(\SL_2(A),X_2(A^2)) \simeq \bigoplus_{\lan a\ran\in \GG_A} H_1(\mu_2(A), \z)
	\simeq \z[\GG_A]\otimes\mu_2(A).
	\]
	Therefore the above exact sequence gives the exact sequence
	\[
	H_1(\SL_2(A), Z_2(A^2)) \arr \z[\GG_A]\otimes\mu_2(A) \arr H_1(\SL_2(A),Z_1(A^2)) 
	\arr \mathcal{RP}_1(A) \arr 0 .
	\]
	
	The exact sequence $0 \arr Z_1(A^2)  \overset{\inc}{\arr} X_1(A^2)  
	\overset{\partial_1}{\arr} X_0(A^2)$ induces the commutative diagram
	\[
	\begin{tikzcd}
		H_1(\SL_2(A), Z_2(A^2))\ar[r]\ar[d] & \z[\GG_A]\otimes\mu_2(A) \ar[r,  "\gamma"]
		\ar[d,"\epsilon\otimes\text{id}_{\mu_2(A)}",twoheadrightarrow] 
		& H_1(\SL_2(A),Z_1(A^2)) \ar[r,"\delta"] \ar[d,"d^1_{2,1}", twoheadrightarrow] & \mathcal{RP}_1(A) \ar[r] \ar[d] 
		& 0\\
		0 \ar[r] & \mu_2(A) \ar[r, equal] & \mu_2(A) \ar[r] & 0 \ar[r] & 0
	\end{tikzcd}
	\]
	By the Snake lemma we have the exact sequence
	\[
	H_1(\SL_2(A), Z_2(A^2)) \arr \II_A\otimes\mu_2(A) \overset{\gamma}{\arr} E^2_{2,1} 
	\arr  \mathcal{RP}_1(A) \arr 0 .
	\]
	
	Let $G$ be a  group and let $g,g'$ be two  commuting elements of $G$. Set
	\[
	\pmb{c}(g,g'):=([g|g']-[g'|g])\otimes 1\in H_2(G,\z)=H_2(B_\bullet(G)\otimes_G \z).
	\]
	
	\begin{lem}\label{gamma-d}
		The composite 
		\[
		\II_A\otimes\mu_2(A) \overset{\gamma}{\larr} E^2_{2,1} \overset{d^2_{2,1}}{\larr}  
		H_2(B(A),\z)\simeq (A^{\times}\wedge A^{\times})\oplus \SS_2
		\]
		sends $\Lan  a\Ran\otimes b$ to 
		$\Big(a\wedge b, \pmb{c}(\mtxx{1}{a+1}{0}{1},\mtxx{b}{0}{0}{b})\Big)$.   
	\end{lem}
	\begin{proof}
		The element $\Lan a\Ran\otimes b\in\II_A\otimes\mu_2(A)$ is represented by 
		$[b]\otimes((\pmb{\infty},\pmb{0},\pmb{a})-(\pmb{\infty},\pmb{0},\pmb{1}))$. 
		Since $\gamma$  is induced by $\partial_2$, we see that 
		$\gamma(\Lan a\Ran\otimes(b))$ is represented by 
		$[b]\otimes\partial_2((\pmb{\infty},\pmb{0},\pmb{a})
		-(\pmb{\infty},\pmb{0},\pmb{1}))\in B_1(\SL_2(A))\otimes Z_1(A^2)$.
		Consider the diagram
		\[
		\begin{tikzcd}
			B_2(\SL_2(A))\otimes X_0(A^ 2) & B_2(\SL_2(A))\otimes X_1(A^2) 
			\ar["\id_{B_2}\otimes \partial_1"',  l] \ar[d, "d_2\otimes \id_{X_1}"] & \\
			& B_1(\SL_2(A))\otimes X_1(A^ 2) & B_1(\SL_2(A))\otimes  
			Z_1(A^2).\ar[l, "\id_{B_1}\otimes \inc"']
		\end{tikzcd}
		\]
		If $X_{a,b}:=[b]\otimes\partial_2((\pmb{\infty},\pmb{0},\pmb{a})
		-(\pmb{\infty},\pmb{0},\pmb{1}))$, then
		\begin{align*}
			(\id_{B_1}\otimes \inc)(X_{a,b})=&[b]\otimes((\pmb{0},\pmb{a})
			-(\pmb{\infty},\pmb{a})-(\pmb{0},\pmb{1})+(\pmb{\infty},\pmb{1}))\\
			=&(g^{-1}_a-h^{-1}_a-g^{-1}_1+h^{-1}_1)[b]\otimes(\pmb{\infty},\pmb{0})\\
			=&(d_2\otimes \id_{X_1})(Z_{a,b}\otimes(\pmb{\infty},\pmb{0}))
		\end{align*}
		where
		\[
		Z_{a,b}:=[g^{-1}_a|b]-[b|g^{-1}_a]-[g^{-1}_1|b]+[b|g^{-1}_1]
		-[h^{-1}_a|b]+[b|h^{-1}_a]+[h^{-1}_1|b]-[b|h^{-1}_1],
		\]
		with $g_z=\mtxx{0}{1}{-1}{z}$ and $h_z=\mtxx{1}{z^{-1}}{0}{1}$ for $z\in\aa$. 
		Applying $\id_{B_2}\otimes \partial_1$ we have 
		\[
		(\id_{B_2}\otimes \partial_1)(Z_{a,b}\otimes(\pmb{\infty},\pmb{0}))
		=(wZ_{a,b}-Z_{a,b})\otimes(\pmb{\infty}).
		\]
		Now $(wZ_{a,b}-Z_{a,b})\otimes 1$ is a representative of $(d^2_{2,1}\circ\gamma)(\Lan a\Ran\otimes b)$.
		We have the following facts:
		\begin{enumerate}
			\item For any $g\in \SL_2(A)$, $h\in B(A)$ and $b,b'\in\mu_2(A)$,
			\[
			\pmb{c}(hg,b)=\pmb{c}(h,b)+\pmb{c}(g,b), \ \ \ \ \ \pmb{c}(h,bb')=\pmb{c}(h,b)+\pmb{c}(h,b').
			\]
			\item For any $g\in \SL_2(A)$, $w([g|b]-[b|g])\otimes 1$ is a representative of 
			$\pmb{c}(wg,b)-\pmb{c}(w,b)$, i.e.
			\[
			\pmb{c}(wg,b)-\pmb{c}(w,b)=\overline{w([g|b]-[b|g])\otimes 1}.
			\]
			\item For any $h\in B(A)$ and $b\in\mu_2(A)$, we have 
			\[
			\pmb{c}(h^{-1},b)=-\pmb{c}(h,b)=\pmb{c}(h,b^{-1})=\pmb{c}(h,b).
			\]
		\end{enumerate}
		Now, for any $z\in\aa$, from the identity $g^{-1}_z=-h_{z^{-1}}w$, we obtain
		\[
		\pmb{c}(g^{-1}_z,b)=\pmb{c}(h_{z^{-1}},b)+\pmb{c}(w,b)+\pmb{c}(-1,b)
		\]
		(by just adding the null element 
		$(d_3\otimes\id)([-h_{a^{-1}}|w|b]+[b|-h_{a^{-1}}|w]-[-h_{a^{-1}}|b|w])$ 
		and using the first fact above). On the other hand, the second fact above gives, 
		for any $z\in\aa$, the equality
		\[
		\overline{w([g^{-1}_z|b]-[b|g^{-1}_z])\otimes 1}=\pmb{c}(wg^{-1}_z,b)-\pmb{c}(w,b).
		\]
		Moreover the formula $wg^{-1}_z=z^{-1}h^{-1}_{z^{-1}}wh^{-1}_z$ and (1) in above gives the equality
		\[
		\overline{w([g^{-1}_z|b]-[b|g^{-1}_z])\otimes 1}=\pmb{c}(z^{-1},b)+\pmb{c}(h^{-1}_{z^{-1}},b)
		+\pmb{c}(wh^{-1}_z,b)-\pmb{c}(w,b).
		\]
		Also using (2) we have
		\[
		\overline{w([h^{-1}_z|b]-[b|h^{-1}_z])\otimes 1}=\pmb{c}(wh^{-1}_z,b)-\pmb{c}(w,b).
		\]
		Now joining all the formulas above we have:
		\begin{align*}
			\overline{(wZ_{a,b}-Z_{a,b})\otimes 1}&=\pmb{c}(a^{-1},b)+\pmb{c}(h^{-1}_{a^{-1}},b)
			-\pmb{c}(h^{-1}_1,b)-\pmb{c}(h_{a^{-1}},b)+\pmb{c}(h_a,b)\\
			&=\pmb{c}(a,b)+\pmb{c}(h_ah_1,b)=\pmb{c}(a,b)
			+\pmb{c}\Big(\!\mtxx{1}{a^{-1}+1}{0}{1},\mtxx{b}{0}{0}{b}\!\Big)
		\end{align*}
		(in the last equality, we use  (1) and (3)). Substituting $a$ with $a^{-1}$ we see that
		\[
		(d^2_{2,1}\circ \gamma)(\Lan a  \Ran \otimes b)
		=\pmb{c}(a,b)+\pmb{c}\Big(\mtxx{1}{a+1}{0}{1},\mtxx{b}{0}{0}{b}\Big).
		\]
	\end{proof}
	
	We believe that the element $\pmb{c}\Big(\!\mtxx{1}{a+1}{0}{1},\mtxx{b}{0}{0}{b}\!\Big)$,
	appearing in the previous lemma, is trivial for many interesting rings. 
	
	For $a\in A$ and $b\in \mu_2(A)$, let 
	$x_a:=\pmb{c}\Big(\!\mtxx{1}{a}{0}{1},\mtxx{b}{0}{0}{b}\!\Big) \in H_2(B(A),\z)$. 
	This element has order $2$ and $x_a=x_{-a}$.
	Since ${\mtxx {c}{0}{0}{c^{-1}}} {\mtxx {1}{a}{0}{1}} {\mtxx {c}{0}{0}{c^{-1}}}^{-1}\!\!\!
	={\mtxx {1} {ac^2} {0} {1}}$, for any $c\in \aa$ we have $x_a=x_{ac^2}$. (In particular 
	$x_{c^2}=x_1$.) Thus 
	\[
	x_{a(c^2-1)}=0, \ \ \ \ x_c=x_{c^{-1}}.
	\]
	For example if $a\in\WW_A$, then $a+1:=\frac{1}{(a-1)}(a^2-1)$ and hence 
	\[
	x_{a+1}=x_{(a-1)^{-1}(a^2-1)}=0.
	\]
	\begin{exa}
		\par (i) If $H_0(\aa, A)=0$, then $A=\lan c^2-1|c\in \aa\ran$. Thus any $a\in A$ is of the 
		form $a=\sum d(c^2-1)$. This implies that $x_a=0$ for any $a\in A$.
		\par (ii) If $2\in \aa$, then for any $a\in \aa$ we have $x_a=x_{2(a/2)}=2x_{(a/2)}=0$.
		\par (iii)  If $F=A$ is a field, then $x_a=0$: If $\char(F)=2$, then $b=1$ and thus $x_a=0$.
		If $\char(F)\neq 2$, then $2\in F^\times$, and the claim follows from (ii).
		\par (iv) If $A$ is a local ring such that $A/\mmm_A$ has at least $3$ elements, then $x_a=0$:
		If $|A/\mmm_A|=3$, then $2 \in \aa$, and thus the claim follows from (ii). If $|A/\mmm_A|>3$,
		then there is $c\in \aa$, such that  $c^2-1\in \aa$. Thus $H_0(\aa, A)=0$ and the claim follows 
		from (i).
		\par (v) Let $A=\z_{(p)}$, where $p$ is a prime. Then $x_{a+1}=0$ for any $a\in \aa$:
		For $p>2$ the claim follows from (iv). Let $p=2$ and let $a=a'/b'\in \z_{(2)}$. Then $a', b'$
		are odd and so $a+1=(a'+b')/b'=2c'$, $c'\in \z_{(2)}$. Now $x_{a+1}=x_{2c'}=2x_{c'}=0$.
		\par (vi) Let $A=\z\pth$, where $p$ is a prime. Then $x_{a+1}=0$ for any $a\in \aa$:
		If $p=2$, then by (ii), $x_{a+1}=x_a=0$. If $p\neq 2$, then  $a=\pm p^n$, $n\in \z$. 
		Now we have $a+1=2c$, where $c\in A$. Thus $x_{a+1}=0$.
		\par (vii) If $\mu_2(A)=1$, then $x_a=0$: Since $\mu_2(A)=1$, we have $b=1$ and thus $x_a=0$.
	\end{exa}
	
	In the rest of this article we will mostly assume that $x_{a+1}=0$ for any $a\in \aa$, i.e.
	\[
	\im (d_{2,1}^2\circ \gamma)=\aa \wedge \mu_2(A).
	\]
	For example in our important results, for technical reasons, we will assume that 
	\[
	H_2(B(A),\z)\simeq H_2(T(A),\z),
	\]
	i.e. $\SS_2=0$. So the above condition will be satisfied.
	
	Now by the above lemma we have the commutative diagram with exact rows
	\begin{equation}\label{diagram1}
		\begin{tikzcd}
			& \II_A\otimes\mu_2(A)  \ar[r,"\gamma"] \ar[d,twoheadrightarrow] & E^2_{2,1} \ar[r,"\delta"] \ar[d] 
			& \mathcal{RP}_1(A) \ar[r] \ar[d] & 0\\
			0 \ar[r] & A^{\times}\wedge\mu_2(A) \ar[r] & (A^{\times}\wedge A^{\times})\oplus \SS_2 \ar[r] & 
			\displaystyle \frac{A^{\times}\wedge A^{\times}}{A^{\times}\wedge\mu_2(A)}  \oplus \SS_2\ar[r] & 0.
		\end{tikzcd}
	\end{equation}
	
	Recall that for any $a\in \aa$, 
	\[
	\psi_1(a):=({\pmb \infty}, {\pmb 0}, {\pmb a})+({\pmb 0}, {\pmb \infty}, {\pmb a})
	-({\pmb \infty}, {\pmb 0}, {\pmb 1})-({\pmb 0}, {\pmb \infty}, {\pmb 1}) \in \RP(A)
	\]
	(see Remark \ref{II-I}). For any $a\in \WW_A$, let 
	\[
	\overline{\psi}_1(a):=[a]+\lan-1\ran[a^{-1}]\in \overline{\RP}(A).
	\]
	Observe that under the natural map $\overline{\RP}(A) \arr \RP(A)$, $\overline{\psi}_1(a)$ maps to $\psi_1(a)$.
	It is easy to check that for any $a\in\WW_A$, the element
	\[
	g(a):=p_{-1}^+[a]+\Lan1-a\Ran \overline{\psi}_1(a)\in \overline{\RP}(A),
	\]
	lies in the subgroup $\overline{\RP}_1(A)$, where $p_{-1}^+=\lan -1\ran +1\in \z[\GG_A]$ (see Remark \ref{II-I}).
	We denote the image of this elements in $\RP_1(A)$ again by $g(a)$. 
	
	\begin {prp}\label{d2}
	Under the composite
	\[
	\RP_1(A) \arr \displaystyle\frac{A^{\times}\wedge A^{\times}}{A^{\times}\wedge\mu_2(A)} \oplus \SS_2 
	\arr \displaystyle\frac{A^{\times}\wedge A^{\times}}{A^{\times}\wedge\mu_2(A)} 
	\]
	we have 
	\[
	g(a)\mt a\wedge (1-a).
	\]
\end{prp}
\begin{proof}
	From the complex $0 \arr Z_1(A^2) \overset{\inc}{\arr} X_1(A^2) \overset{\partial_1}{\arr} X_0(A^2)  
	\arr 0$ we obtain the first quadrant spectral sequence
	\[
	\EE^1_{p.q}=\left\{\begin{array}{ll}
		H_q(\GL_2(A),X_p(A^2)) & p=0,1\\
		H_q(\GL_2(A),Z_1(A^2)) & p=2\\
		0 & p>2
	\end{array}
	\right.
	\Longrightarrow H_{p+q}(\GL_2(A),\z).
	\]
	This spectral sequence have been studied in \cite[\S 3]{mirzaii2011}. Let 
	$\PP(A):=H_0(\GL_2(A), Z_2(A^2))$. We have a $\z[\GG_A]$-map $\RP(A) \arr \PP(A)$, 
	where $\PP(A)$ has the trivial action of $\GG_A$. Under this map $g(a) \mt 2[a]$. 
	This induces a map  $\RP_1(A) \arr \PP(A)$. One can show that $\EE_{2,1}^2\simeq \PP(A)$ 
	(see \cite[Lemma~3.2]{mirzaii2011}). The map $\SL_2(A) \arr \GL_2(A)$ induces the morphism 
	of spectral sequences
	\[
	\begin{tikzcd}
		E^1_{p,q} \ar[d] \ar[r,Rightarrow] & H_{p+q}(\SL_2(A),\z)\ar[d]\\
		\EE^1_{p,q} \ar[r,Rightarrow] & H_{p+q}(\GL_2(A),\z).
	\end{tikzcd}
	\]
	From this we obtain the commutative diagram
	\[
	\begin{tikzcd}
		E_{2,1}^2 \ar[d] \ar[r, "d_{2,1}^2 "] \ar[d]& H_2(B(A),\z) \ar[d] \ar[r]  &H_2(T(A),\z)\ar[d]\\
		\PP(A)\ar[r,"d_{2,1}^2 "]&\displaystyle\frac{H_2(B_2(A),\z)}{d_{1,2}^2(H_2(T_2(A),\z))}\ar[r] 
		&\displaystyle\frac{H_2(T_2(A),\z)}{d_{1,2}^2(H_2(T_2(A),\z))},
	\end{tikzcd}
	\]
	where
	\[
	B_2(A):=\stabe_{\GL_2(A)}(\pmb{\infty})=\Bigg\{\begin{pmatrix}
		a & b\\
		0 & d
	\end{pmatrix}:a, d\in \aa, b\in A\bigg\},
	\]
	\[
	T_2(A):=\stabe_{\GL_2(A)}({\pmb \infty},{\pmb 0})=\Bigg\{\begin{pmatrix}
		a & 0\\
		0 & d
	\end{pmatrix}:a, d\in \aa\bigg\}.
	\]
	This together with diagram (\ref{diagram1}) induce the commutative diagram
	\[
	\begin{tikzcd}
		\RP_1(A) \ar[d] \ar[r] \ar[d]&\displaystyle\frac{\aa\wedge \aa}{\aa\wedge \mu_2(A)} \ar[d, hook]\\
		\PP(A) \ar[r]  & (\aa\wedge\aa) \oplus S_\z^2(A),
	\end{tikzcd}
	\]
	where $S_\z^2(A):=(\aa \otimes \aa)/\lan a\otimes b+b\otimes a:a,b\in \aa \ran$. The vertical 
	map on the right is given by $a\wedge b\arr (2a\wedge b, 2(a\otimes b))$ and the bottom horizontal map 
	is given by $[a]\mt (a\wedge (1-a), -a \otimes (1-a))$ \cite[Lemma 4.1]{mirzaii2011}. Now the claim follows 
	from the fact that the composite
	\[
	\RP_1(A)  \arr \PP(A) \arr  (\aa\wedge\aa) \oplus S_\z^2(A)
	\]
	maps $g(a)$ to $2(a\wedge (1-a), -a \otimes (1-a))$.
\end{proof}

We denote the differential $d_{2,1}^2$ by $\lambda_1$:
\[
\lambda_1: \RP_1(A) \arr H_2(B(A),\z)\simeq \displaystyle
\frac{A^{\times}\wedge A^{\times}} {A^{\times}\wedge\mu_2(A)}\oplus \SS_2 .
\]
The kernel of $\lambda_1$ is called the {\it refined Bloch group} of $A$ 
and is denoted by $\RB(A)$:
\[
\RB(A):=\ker(\lambda_1).
\]
From the spectral sequence we obtain a natural surjective map 
\[
H_3(\SL_2(A),\z) \two \RB(A).
\]
Let $\Sigma_2'=\{1, \sigma'\}$ be the symmetric group of order $2$. This group acts on $\tors(\mu(A),\mu(A))$ as 
$(\sigma', x)\mapsto -\sigma_1(x)$, where $\sigma_1:\tors(\mu(A),\mu(A))\arr \tors(\mu(A),\mu(A))$ is 
obtained by interchanging the group $\mu(A)$.
\begin{thm}[Refined Bloch-Wigner in $\char=2$ \cite{B-E-2023}]\label{B-E2023}
	Let $A$ be a ring such that 
	\par {\rm (i)} $\mu_2(A)=1$,
	\par {\rm (ii)} $X_\bullet(A^2) \arr \z$ is exact in dimension $< 2$
	\par {\rm (iii)} $H_3(T(A), \z) \simeq H_3(B(A), \z)$.\\
	Then we have the exact sequence
	\[
	\tors(\mu(A), \mu(A))^{\Sigma_2'} \arr H_3(\SL_2(A),\z) \arr \RB(A) \arr 0.
	\]
	If $A$ is a domain, then we have the exact sequence
	\[
	0\arr \tors(\mu(A), \mu(A)) \arr H_3(\SL_2(A),\z) \arr \RB(A) \arr 0.
	\]
\end{thm}
\begin{proof}
	This is a slight generalization of \cite[Theorem 6.1]{B-E-2023} and the proof is the same.
\end{proof}

We now study the map $\II_A\otimes\mu_2(A)\arr A^{\times}\wedge\mu_2(A)\se \aa \wedge \aa$ given by
$\Lan a\Ran \otimes b\mt a\wedge b$ (when $A$ is a domain). Clearly $\II^2_A\otimes\mu_2(A)$ 
is in the kernel of this map. This induces the map
\[
\GG_A\otimes \mu_2(A)\simeq (\II_A/\II_A^2) \otimes\mu_2(A)\arr A^{\times}\wedge\mu_2(A), 
\]
\[
\lan a\ran\otimes b \mt \Lan a \Ran \otimes b \mt  a\wedge b.
\]

\begin{lem}\label{z2'}
	Let $A$ be a domain. Then the kernel of the map $\GG_A \otimes \mu_2(A)\arr A^{\times}\wedge\aa$, 
	given by $\lan a\ran \otimes (-1) \mt a\wedge(-1)$, has at most two elements.
\end{lem}
\begin{proof}
	We may assume that $\char(A)\neq 2$. In this case $\GG_A \otimes \mu_2(A)\simeq \GG_A$.
	Let $a\wedge(-1)=0$ in $A^{\times}\wedge A^{\times}$. We know 
	that $A^{\times}=\varinjlim H$, where $H$ runs through all finitely generated subgroups 
	of $A^{\times}$. As the direct limit commutes with wedge product, we have
	$A^{\times}\wedge A^{\times}=\varinjlim H\wedge H$. We may take a finitely 
	generated subgroup $H$ such  $a,-1\in H$ and $a\wedge (-1)=0\in H\wedge H$.
	
	Let $H\simeq F \times T$, where $F$ is torsion free and $T$ is a finite cyclic group. 
	Thus $-1\in T$ and we have
	\[
	H\wedge H \simeq (F\wedge F)\oplus (F\otimes T)\oplus (T\wedge T).
	\] 
	Clearly $T\wedge T=0$.
	Let $a=p\omega$ with $p\in F$ and $\omega \in T$. From $a\wedge(-1)=0\in H \wedge H$, 
	it follows that $p\otimes(-1)=0$ and $\omega\wedge(-1)=0$. As $-1\in T$, $T$ has even order.
	Thus $p\otimes(-1)=0$ implies that $p$ is a square. Therefore $\lan a\ran=\lan\omega\ran$. 
	This completes the proof.
\end{proof}

Now let $A$ be a domain. Then from the  commutative diagram (\ref{diagram1}), we obtain the exact sequence
\[
H_1(\SL_2(A), Z_2(A^2)) \arr J \overset{\gamma}{\arr} E^3_{2,1} \arr  \RB(A) \arr 0,
\]
where $J$ sits in the exact sequence $\II_A^2\otimes\mu_2(A)\arr J \arr (\z/2)' \arr 0$
with $(\z/2)'$ a subgroup of $\z/2$ (Lemma \ref{z2'}).

\section{The low dimensional homology of \texorpdfstring{$\SM_2$}{Lg} }

Let $\SM_2(A)$ denotes the group of monomial matrices in $\SL_2(A)$. So $\SM_2(A)$ consists of matrices 
$\begin{pmatrix}
	a & 0\\
	0 & a^{-1}
\end{pmatrix}$
and 
$\begin{pmatrix}
	0 & a\\
	-a^{-1} & 0
\end{pmatrix}$,
where $a\in \aa$. Let $\hat{X}_0(A^2)$ and $\hat{X}_1(A^2)$ be the free $\z$-modules generated by the sets 
\[
\SM_2(A)({\pmb \infty}):=\{g.({\pmb \infty}):g\in \SM_2(A)\}, \ \ \ \ 
\SM_2(A)(\pmb{\infty}, {\pmb 0}):=\{g.({\pmb \infty}, {\pmb 0}):g\in \SM_2(A)\},
\]
respectively. It is easy to see that the sequence of $\SM_2(A)$-modules
\[
\hat{X}_1(A^2)  \overset{\hat{\partial}_1}{\arr} \hat{X}_0(A^2) \overset{\hat{\epsilon}{}}{\arr} \z \arr 0
\]
is exact and
\[
\ker(\hat{\partial}_1)=\z\{(\pmb{\infty},\pmb{0})+(\pmb{0}, \pmb{\infty})\}.
\]
We denote this kernel by $\hat{Z}_1(A^2)$. Observe that $\hat{Z}_1(A^2)\simeq \z$ 
and $\SM_2(A)$ acts trivially on it. From the complex
\begin{equation}\label{comp2}
	0 \arr\hat{Z}_1(A^2)\overset{\hat{\inc}}{\arr}\hat{X}_1(A^2)\overset{\hat{\partial}_1}{\arr}\hat{X}_0(A^2)\arr 0,
\end{equation}
we obtain the first quadrant spectral sequence
\[
\hat{E}^1_{p.q}=\left\{\begin{array}{ll}
	H_q(\SM_2(A),\hat{X}_p(A^2)) & p=0,1\\
	H_q(\SM_2(A),\hat{Z}_1(A^2)) & p=2\\
	0 & p>2
\end{array}
\right.
\Rightarrow H_{p+q}(\SM_2(A),\z).
\]
Since the  complex (\ref{comp2}) is a $\SM_2(A)$-subcomplex of (\ref{comp1}), we have a natural 
morphism of spectral sequences
\begin{equation}\label{seqspmorph}
	\begin{tikzcd}
		\hat{E}^1_{p,q} \ar[d] \ar[r,Rightarrow] & H_{p+q}(\SM_2(A),\z)\ar[d]\\
		E^1_{p,q} \ar[r,Rightarrow] & H_{p+q}(\SL_2(A),\z).
	\end{tikzcd}
\end{equation}
As in case of $\SL_2(A)$, we have $\hat{X}_0\simeq \Ind _{T(A)}^{\SM_2(A)}\z$ and 
$\hat{X}_1\simeq \Ind _{T(A)}^{\SL_2(A)}\z$. Thus  by Shapiro's lemma we have
\[
\hat{E}_{0,q}^1 \simeq H_q(T(A),\z), \ \ \ \ \
\hat{E}_{1,q}^1 \simeq H_q(T(A),\z).
\]
Therefore
\[
\hat{E}^1_{p.q}=\left\{\begin{array}{ll}
	H_q(T(A),\z) & p=0,1\\
	H_q(\SM_2(A),\z) & p=2\\
	0 & p>2
\end{array}
\right.
\Rightarrow H_{p+q}(\SM_2(A),\z).
\]
Moreover, $\hat{d}_{1, q}^1=H_q(\hat{\sigma}) - H_q(\hat{\inc})=\hat{\sigma}_\ast-\hat{\inc}_\ast$, where
$\hat{\sigma}: T(A) \arr T(A)$ is given by $X \arr wXw^{-1}=X^{-1}$. Thus  $\hat{d}_{1,0}^1$ is trivial, 
$\hat{d}_{1,1}^1$ is induced by the map $X\mt X^{-2}$ and $\hat{d}_{1,2}^1$ is trivial. 

A direct calculation shows that the map $\hat{d}_{2,q}: H_q(\SM_2(A),\z) \arr H_q(T(A),\z)$ is the 
transfer map \cite[\S 9, Chap. III]{brown1994}. Hence the composite
\[
H_q(\SM_2(A),\z) \overset{\hat{d}_{2,q}}{\arr} H_q(T(A),\z)  \overset{\inc_\ast}{\arr} H_q(\SM_2(A),\z)
\]
coincides with multiplication by 2 \cite[Proposition 9.5, Chap. III]{brown1994}. 
In particular, $\hat{d}_{2,0}:\z \arr \z $ is multiplication by 2. 
From these we obtain the exact sequence 
\[
1 \arr \GG_A \arr H_1(\SM_2(A),\z)\arr \z/2 \arr 0. 
\]
If fact this can be obtain directly from the extension 
$1 \arr T(A) \arr \SM_2(A) \arr \lan \overline{w}\ran \arr 1$: 
\[
1 \arr \GG_A \arr H_1(\SM_2(A),\z) \arr \lan \overline{w}\ran \arr 1.
\]
Observe that $w^2=\begin{pmatrix}
	-1&0\\
	0&-1
\end{pmatrix} \in T(A)$.  
A direct calculation shows that $\hat{d}_{2,1}^1(\overline{w})=-1$ and 
$\hat{d}_{2,1}^1\mid_{\GG_A}=0$. Thus 
\[
\hat{E}_{1,1}^2=\mu_2(A)/\{\pm 1\}, \ \ \ \ \ \hat{E}_{2,1}^2=\GG_A.
\]
Again a direct calculation shows that 
\[
\hat{d}_{2,1}^ 2:\GG_A \arr H_2(T(A),\z)\simeq \aa \wedge \aa
\]
is given by $\lan a\ran \mapsto a\wedge (-1)$. Therefore from the spectral sequence 
$\hat{E}_{p,q}^ 1 \Rightarrow H_{p+q}(\SM_2(A),\z)$ we obtain the exact sequence 
\[
0 \arr \displaystyle\frac{\aa \wedge \aa}{\aa \wedge \{\pm 1\}} \arr H_2(\SM_2(A),\z) 
\arr \mu_2(A)/\{\pm 1\} \arr 1.
\]
Thus we have:

\begin{lem}\label{H2SM2}
	If $\mu_2(A)=\{\pm 1\}$, then 
	$H_2(\SM_2(A),\z)\simeq \displaystyle\frac{\aa \wedge \aa}{\aa \wedge \mu_2(A)}$.
\end{lem}

Now if $\mu_2(A)=\{\pm 1\}$, then it follows from this lemma that the image of the map
$\hat{d}_{2,2}^1:H_2(\SM_2(A),\z)  \arr \aa\wedge \aa$ is $2(\aa\wedge\aa)$. Thus 
$\hat{E}_{1,2}^2\simeq\displaystyle\frac{\aa\wedge\aa}{2(\aa\wedge\aa)}$.
Moreover one can show that 
$\hat{E}_{2,2}^2\simeq\displaystyle \frac{{}_2(\aa \wedge \aa)}{\aa\wedge \mu_2(A)}$.

\section{The third homology of \texorpdfstring{$\SL_2$}{Lg}}

Assume that $A$ satisfies the condition that $X_\bullet(A^2)\arr \z$ is exact in dimension $<2$.
Then the natural map $\alpha: \GG_A=\hat{E}_{2,1}^2 \arr E_{2,1}^2$ sits in the diagram
\[
\begin{tikzcd}
	& \GG_A \ar[d, "\alpha"] & \\
	\II_A\otimes \mu_2(A) \ar[r, "\gamma"] & E_{2,1}^2 \ar[r, "\delta"] 
	& \RP_1(A) \larr 0.
\end{tikzcd}
\]


\begin{lem}\label{image}
	The composite map $\delta\circ\alpha:\GG_A \arr \RP_1(A)$ is given by $\lan a\ran \mapsto \psi_1(a^2)$.
\end{lem}
\begin{proof}
	The element $\lan a\ran\in\GG_A$ is represented by 
	\[
	[a]\otimes\{({\pmb \infty},{\pmb 0})+({\pmb 0},{\pmb \infty})\}\in H_1(\SM_2(A),\hat{Z}_1(A^2)).
	\]
	Its image in $H_1(\SL_2(A),Z_1(A^2))$, through $\alpha$, is represented by the element
	\[
	S:=[a]\otimes\partial_2((\pmb{\infty},{\pmb 0},{\pmb a^2})+({\pmb 0},{\pmb \infty}, {\pmb a^2})).
	\]
	We have
	\begin{align*}
		\delta(S) & =(d_1\otimes\text{id}_{Z_2(X^2)})\bigg ([a]\otimes\partial_2(({\pmb \infty},
		{\pmb 0}, {\pmb a^2})+({\pmb 0},{\pmb \infty}, {\pmb a^2}))\bigg)\\
		&=[\ ]\otimes\bigg(({\pmb \infty},{\pmb 0},{\pmb 1})+({\pmb 0},{\pmb \infty},{\pmb 1})-({\pmb \infty},{\pmb 0},
		{\pmb a^2})-({\pmb 0},{\pmb \infty},{\pmb a^2})\bigg.
	\end{align*}
	It is straightforward to check that this element represents $-\psi_1(a^2)$. Thus
	\[
	\delta(S)=-\psi_1(a^2)=\psi_1(a^2)
	\]
	(see \cite[Lemma 3.21(5)]{C-H2022}).
\end{proof}

For any $a\in A^{\times}$, let  $X_a$ and $X'_a$ denote the elements 
$({\pmb \infty},{\pmb 0}, {\pmb a})$ and $({\pmb 0},{\pmb \infty}, {\pmb a})$ of 
$X_2(A^2)$, respectively. Let $\chi_a \in H_1(\SL_2(A), Z_1(A^2))$ be 
represented by $[wa]\otimes\partial_2(X_{-a}-X_a)$, where $w=E(0)={\mtxx{0}{1}{-1}{0}}$.  We usually write
\[
\chi_a:= [wa]\otimes\partial_2(X_{-a}-X_a).
\]
Recall that for simplicity 
$\begin{pmatrix}
	a & 0\\
	0 & a^{-1}
\end{pmatrix}$ usually is denoted by $a$.

\begin{lem}\label{crucial0}
	For any $a\in \aa$, $\gamma(\Lan a\Ran\otimes(-1))-\alpha(\langle a\rangle)=\lan -1\ran\Lan a\Ran.\chi_1$.
\end{lem}
\begin{proof}
	Let $Y:=({\pmb \infty},{\pmb 0})+({\pmb 0},{\pmb \infty})\in Z_1(A^2)$. 
	For any $a\in A^{\times}$, we have
	\begin{itemize}
		\item [(a)] $d_2([wa|wa])=wa[wa]-[-1]+[wa]$,
		\item[(b)] $d_2([w|a])=w[a]-[wa]+[w]$.
	\end{itemize}
	Thus modulo $\text{im}(d_2\otimes\text{id}_{Z_1(A^2)})$, we  have
	\begin{enumerate}
		\item $[-1]\otimes\partial_2(X_{-a})=[wa]\otimes\partial_2(X'_a)+[wa]\otimes\partial_2(X_{-a})$,
		\item $[wa]\otimes Y=[a]\otimes Y+[w]\otimes Y$.
	\end{enumerate}
	Hence
	\begin{align*}
		[wa]\otimes\partial_2(X_{-a}-X_a)=&[wa]\otimes\partial_2(X_{-a})-[wa]\otimes\partial_2(X_a)\\
		=&[-1]\otimes\partial_2(X_{-a})-[wa]\otimes\partial_2(X'_a)-[wa]\otimes\partial_2(X_a)\\
		=&[-1]\otimes\partial_2(X_{-a})-[wa]\otimes Y\\
		=&[-1]\otimes\partial_2(X_{-a})-([a]\otimes Y+[w]\otimes Y)\\
		=&[-1]\otimes\partial_2(X_{-a})-[w]\otimes Y-\alpha(\lan a\ran)\\
		=&[-1]\otimes\partial_2(X_{-a})-[w]\otimes\partial_2(X_1+X'_1)-\alpha(\lan a\ran).
	\end{align*}
	Now, using  the identity (1)  in above for $a=1$, we get
	\begin{align*}
		[wa]\otimes\partial_2(X_{-a}-X_a)-[w]\otimes\partial_2(X_{-1}-X_1)&
		=[-1]\otimes\partial_2(X_{-a}-X_{-1})-\alpha(\lan a\ran)\\
		&=\lan-1\ran\gamma(\Lan a\Ran\otimes(-1))-\alpha(\lan a\ran).
	\end{align*}
	On the other hand,
	\begin{align*}
		[wa]\otimes\partial_2(X_{-a}\!\!-X_a)\!-\![w]\otimes\partial_2(X_{-1}\!-\!X_1)
		&= \lan a\ran([w]\otimes\partial_2(X_{-1}\!-\!X_1))-[w]\otimes\partial_2(X_{-1}\!-\!X_1)\\
		&=  \Lan a\Ran([w]\otimes\partial_2(X_{-1}-X_1))\\
		&=\Lan a\Ran \chi_1.
	\end{align*}
	Therefore $\Lan a\Ran\cdot\chi_1=\lan-1\ran\gamma(\Lan a\Ran\otimes(-1))-\alpha(\lan a\ran)$.
\end{proof}

\begin{rem}
	It is straightforward to show that 
	\[
	\delta(\chi_1)=\psi_1(-1)\in \RP_1(A). 
	\]
\end{rem}

\begin{cor}\label{crucial}
	If $-1\in (A^{\times})^2$, then for any $a\in A^{\times}$, 
	$\gamma(\Lan a\Ran\otimes(-1))=\alpha(\langle a\rangle)$.
\end{cor}
\begin{proof}
	First observe that for any $s\in \aa$ and $X\in X_2(A^2)$, we have
	\[
	[w]\otimes(sX-X)=[s]\otimes(wX+sX).
	\]
	Now if $i^2=-1$, then by the above relation we have
	\begin{align*}
		[w]\otimes\partial_2(X_{-1}-X_1) &=[w]\otimes\partial_2(iX_1-X_1)\\
		& =[i]\otimes\partial_2(wX_1+iX_1)\\
		& =[i]\otimes\partial_2(X'_1+X_1)\\
		& =[i]\otimes Y=\alpha (\lan i\ran).
	\end{align*}
	Now the claim follows from Lemma \ref{crucial0}.
\end{proof}

\begin{cor}\label{exact2}
	Let $\mu_2(A)=\{\pm 1\}$ and $-1\in (A^{\times})^2$. Then $\gamma(\II_A^2\otimes \mu_2(A))=0$. 
	In particular, we have the exact sequence
	\[
	\GG_A \overset{\alpha}{\larr} E_{2,1}^2 \overset{\delta}{\larr} \RP_1(A) \arr 0.
	\]
\end{cor}
\begin{proof}
	The ideal $\II_A^2$ is generated by the elements $\Lan a\Ran\Lan b\Ran=\Lan ab\Ran-\Lan a\Ran-\Lan b\Ran$. 
	Thus by the above corollary
	\[
	\gamma(\Lan a\Ran\Lan b\Ran\otimes (-1))=\alpha(\lan ab\ran)-\alpha(\lan a\ran)-\alpha(\lan b\ran)
	=\alpha( \lan aba^{-1}b^{-1}\ran=\alpha(\lan 1\ran)=0.
	\]
	The second part follows from the first part and the fact that $\II_A/\II_A^2\simeq \GG_A$
	and $\im(\gamma)=\im(\alpha)$.
\end{proof}

\begin{thm}\label{prthm}
	Let $A$  be a commutative ring such that
	\par {\rm (i)} $\mu_2(A)=\{\pm 1\}$ and $-1\in (A^{\times})^2$,
	\par {\rm (ii)} $X_\bullet(A^2) \arr \z$ is exact in dimension $<2$.
	\par {\rm (iii)} $H_i(T(A),\z)\simeq H_i(B(A),\z)$ for $i=2,3$. \\
	Then we have the exact sequence
	\[
	H_3(\SM_2(A),\mathbb{Z}) \arr H_3(\SL_2(A),\mathbb{Z}) \arr \RB(A) \arr 0.
	\]
\end{thm}
\begin{proof}
	The morphism of spectral sequences (\ref{seqspmorph}) induces a map of filtration
	\[
	\begin{array}{cccccccc}
		0\se & \hat{F}_0 & \se &\hat{F}_1   &\se &\hat{F}_2   &\se & \hat{F}_3=H_3(\SM_2(A),\z)\\
		&\downarrow &      &\downarrow &    & \downarrow &    & \downarrow\\
		0\se &    F_0    &\se   &  F_1      &\se & F_2        &\se & F_3 =H_3(\SL_2(A),\z)
	\end{array}
	\]
	where $E_{p,3-p}^\infty=F_p/F_{p-1}$  and $\hat{E}_{p,3-p}^\infty=\hat{F}_p/\hat{F}_{p-1}$.
	Clearly $F_2=F_3$ and $\hat{F}_2=\hat{F}_3$. Consider the following commutative 
	diagram with exact rows 
	\begin{equation}\label{diagram}
		\begin{tikzcd}
			0 \ar[r] & \hat{F}_1 \ar[r] \ar[d] & H_3(\SM_2(A),\z) \ar[r] \ar[d,"\inc_\ast "] 
			& \hat{E}^{\infty}_{2,1} \ar[r] \ar[d] & 0 \\
			0 \ar[r] & F_1 \ar[r] & H_3(\SL_2(A),\z) \ar[r] & E^{\infty}_{2,1} \ar[r] & 0.
		\end{tikzcd}
	\end{equation}
	By Corollary \ref{exact2}, we have the exact sequence 
	$\hat{E}^2_{2,1}\arr E^2_{2,1} \arr \RP_1(A) \arr 0$.
	From the commutative diagram with exact rows
	\[
	\begin{tikzcd}
		& \hat{E}^2_{2,1}\ar[r]\ar[d,twoheadrightarrow,"\hat{d}_{2,1}^2"] & E^2_{2,1}\ar[r] 
		\ar[d, "d_{2,1}^2"] & \mathcal{RP}_1(A) \ar[r] \ar[d] & 0\\
		0 \ar[r] & A^{\times}\wedge\mu_2(A) \ar[r] & (A^{\times}\wedge A^{\times}) \ar[r] & 
		\displaystyle\frac{A^{\times}\wedge A^{\times}}{A^{\times}\wedge\mu_2(A)} \ar[r] & 0
	\end{tikzcd}
	\]
	we obtain the exact sequence
	\[
	\hat{E}^{\infty}_{2,1}\arr E^{\infty}_{2,1} \arr \mathcal{RB}(A) \arr 0.
	\]
	Now consider the commutative diagram with exact rows
	\[
	\begin{tikzcd}
		0 \ar[r] & \hat{F}_0 \ar[r] \ar[d] & \hat{F}_1 \ar[r] \ar[d] & \hat{E}^{\infty}_{1,2} \ar[r] \ar[d] & 0 \\
		0 \ar[r] & F_0 \ar[r] & F_1 \ar[r] & E^{\infty}_{1,2} \ar[r] & 0.
	\end{tikzcd}
	\]
	Since $\hat{E}^1_{0,3}\simeq E^1_{0,3} $, the natural map  $\hat{F}_0 \arr F_0$ is surjective. Moreover, 
	since $\hat{E}^1_{1,2}\simeq E^1_{1,2} $, the map $\hat{E}^{\infty}_{1,2}\arr E^{\infty}_{1,2}$ is 
	surjective. These imply that the map $\hat{F}_1\arr F_1$ is surjective. Now the claim follows by 
	applying the snake lemma to the diagram (\ref{diagram}).
\end{proof}

\begin{rem}\label{-1-cond}
	We think that the condition $-1\in \aa^2$ in Theorem~\ref{prthm} is not essential (at least when $A$ is a domain). 
	To remove this condition we need to prove that under the map $\gamma: \II_A\otimes \mu_2(A) \arr E_{2,1}^2$, 
	$\II_A^2 \otimes \mu_2(A)$ maps to zero. Having this, then 
	\[
	\GG_A\simeq \GG_A\otimes \mu_2(A) \overset{\bar{\gamma}}{\larr} E_{2,1}^2 \arr \aa \wedge \mu_2(A)
	\ \ \text{and} \ \ \GG_A \overset{\alpha}{\larr} E_{2,1}^2 \arr \aa \wedge \mu_2(A)
	\]
	have the same kernel by Lemma~\ref{z2'}. Then we can proceed as in the above proof.
\end{rem}

\begin{exa}
	Here we give examples of rings  that satisfy the conditions of Theorem~\ref{prthm}:
	\par (1) 
	Any local domain of characteristic $2$ such that its residue field has more than $64$
	elements (Proposition \ref{GE2C}, Theorem \ref{iso-hut}).
	\par (2) 
	Let $B$ be a domain such that $-1$ is square. Let $\ppp$ be a prime ideal of $B$ such that 
	either $B/\ppp$ is infinite or if $|B/\ppp|=p^d$, then $(p-1)d>6$. Then $A:=B_\ppp$ satisfies 
	in the conditions of Theorem ~\ref{prthm} (see Proposition \ref{GE2C} and Theorem \ref{iso-hut}).
	\par (3) 
	Any domain with many units such that $-1$ is an square (e.g $F$-algebras which are domains and $F$ 
	is an algebraically closed) \cite[\S 2]{mirzaii2011}.
\end{exa}
\section{A spectral sequence for relative homology}

Let $G$ be a group and $M$ a $G$-module. We denote these by a pair $(G,M)$. A morphism of pairs
$(f, \sigma): (G', M') \arr (G, M)$ is a pair of group homomorphisms $f:G' \arr G$ and $\sigma: M'\arr M$ 
such that
\[
\sigma(g'm')=f(g')\sigma(m').
\]
This means that $\sigma$ is a map of $G'$-modules.

For a group $H$ let $C_\bullet(H) \arr \z$ be the standard resolution of $\z$ over $\z[H]$ \cite[Chap.I, \S 5]{brown1994}. 
The map $f:G'\arr G$, induces in a natural way a morphism of complexes 
$f_\bullet:C_\bullet(G')\arr C_\bullet(G)$.

The morphism of the pairs $(f, \sigma):(G', M') \arr (G, M)$, induces a morphism of complexes
\[
f_\bullet\otimes \sigma: C_\bullet(G')\otimes_{G'}M' \to C_\bullet(G)\otimes_{G}M.
\]

Let $G'$ be a subgroup of $G$ and $M'$ be a $G'$-submodule of $M$. We take $(i,\sigma): (G',M')\harr (G,M)$ as
the natural pair of inclusion maps. Then the morphism
\[
i_\bullet\otimes \sigma: C_\bullet(G')\otimes_{G'}M' \to C_\bullet(G)\otimes_{G}M
\]
is injective. We denote the $n$-homology of the quotient complex 
$C_\bullet(G)\otimes_{G}M/C_\bullet(G')\otimes_{G'}M'$ by $H_n(G, G'; M', M)$:
\[
H_n(G, G'; M, M'):=H_n(C_\bullet(G)\otimes_{G}M/C_\bullet(G')\otimes_{G'}M').
\]
If $M'=M$, then $H_n(G, G'; M, M')$ is the usual relative homology group $H_n(G, G'; M)$.

From the exact sequence of complexes
\[
0 \arr C_\bullet(G')\otimes_{G'}M' \to C_\bullet(G)\otimes_{G}M \arr 
C_\bullet(G)\otimes_{G}M/C_\bullet(G')\otimes_{G'}M'\arr 0
\]
we obtain the long exact sequence
\[
\cdots \arr H_n(G', M') \arr H_n(G, M) \arr H_n(G, G'; M, M') \arr H_{n-1}(G', M') 
\]
\[
\arr H_{n-1}(G, M)\arr H_{n-1}(G, G'; M, M') \arr \cdots
\]

\begin{prp}\label{relative}
	Let $G'$ be a subgroup of $G$. Let $L_\bullet' \arr M'$ be an exact $G'$-subcomplex of an exact $G$-complex 
	$L_\bullet \arr M$. Then we have the first quadrant spectral sequence
	\[
	\Eb_{p,q}^1= H_q(G,G'; L_p, L_{p}') \Rightarrow H_{p+q}(G, G'; M, M').
	\]
\end{prp}
\begin{proof}
	Let $i:G'\harr G$ and $\sigma_\bullet:L_\bullet' \harr L_\bullet$ be the usual inclusions. 
	The morphism of double complexes
	\[
	i_\bullet\otimes\sigma_\bullet:C_\bullet(G')\otimes_{G'} L_\bullet' \arr C_\bullet(G)\otimes_G L_\bullet
	\]
	is injective. We denote its quotient by $D_{\bullet,\bullet}$: 
	$D_{\bullet,\bullet}=\coker(i_\bullet\otimes\sigma_\bullet)$.
	These double complexes induces two spectral sequences
	\[
	\EE_{p,q}^1 (I) = H_q(D_{p,\bullet}) \Rightarrow H_{p+q}(\Tot(D_{\bullet,\bullet})),
	\ \ \ \ 
	\EE_{p,q}^1 (II) = H_q(D_{\bullet,p}) 
	\Rightarrow H_{p+q}(\Tot(D_{\bullet,\bullet})).
	\]
	These are the spectral sequences
	\[
	\EE_{p,q}^1 (I) = H_q\Bigg(\frac{C_p(G)\otimes_G L_\bullet}{C_p(G')\otimes_{G'} L_\bullet'}\Bigg) 
	\Rightarrow H_{p+q}(\Tot(D_{\bullet,\bullet})),
	\]
	and
	\[
	\EE_{p,q}^1(II) = H_q\Bigg(\frac{C_\bullet \otimes_G L_p}{C_\bullet' \otimes_{G'} L_p'}\Bigg) 
	\Rightarrow H_{p+q}(\Tot(D_{\bullet,\bullet})).
	\]
	By definition $\EE_{p,q}^1(II)=H_q(G,G', L_p,L_p')$. Moreover since $L_\bullet$ and $L_\bullet'$ are exact
	in dimension $>0$, we have $\EE_{p,q}^1(I)=0$ for any $q>0$. For $q=0$, we have
	$\EE_{p,0}^1(I)\simeq \displaystyle\frac{C_p(G)\otimes_G M}{C_p(G')\otimes_{G'} M'}$.
	The homology of the sequence $\EE_{p+1,0}^1(I)\arr \EE_{p,0}^1(I)\arr \EE_{p,0}^1(I)$ is
	\[
	\EE_{p,0}^2(I)\simeq H_q(G, G', M, M').
	\]
	Now by an easy analysis of the spectral sequence $\EE_{p,q}^1(I)$, for any $n\geq 0$ we obtain the isomorphism
	\[
	H_n(\Tot(D_{\bullet,\bullet}))\simeq H_n(G, G'; M, M').
	\]
	Thus if we take $\Eb_{p,q}^1:=\EE_{p,q}^1(II)$, then we obtain the spectral sequence
	\[
	\Eb_{p,q}^1=H_q(G,G'; L_p,L_p') \Rightarrow H_{p+q}(G, G'; M, M').
	\]
\end{proof}

\section{The groups \texorpdfstring{$\RP_1(A)$}{Lg} and  \texorpdfstring{$H_3(\SL_2(A),\SM_2(A),\z)$}{Lg}}\label{RP-SL}

Let $\mu_2(A)=\{\pm1\}$. Assume that $A$ satisfies the condition that $X_\bullet(A^2)\arr \z$ is
exact in dimension $<1$. The complex 
\[
0\arr\hat{Z}_1(A^2) \arr \hat{X}_1(A^2)\arr\hat{X}_0(A^2)\arr 0
\]
is a $\SM_2(A)$-subcomplex of the $\SL_2(A)$-complex 
\[
0 \arr Z_1(A^2) \arr X_1(A^2)  \arr X_0(A^2) \arr 0.
\]
By Proposition \ref{relative}, from the morphism of complexes
\[
\begin{tikzcd}
	0\ar[r]&\hat{Z}_1(A^2) \ar[r]\ar[d]&\hat{X}_1(A^2)\ar[r]\ar[d]&\hat{X}_0(A^2)\ar[r]\ar[d] &0\\
	0 \ar[r] & Z_1(A^2)  \ar[r]  & X_1(A^2)  \ar[r]& X_0(A^2) \ar[r] & 0,
\end{tikzcd}
\]
we obtain the first quadrant spectral sequence
\[
\Eb_{p,q}^ 1\!=\!
\begin{cases}
	H_q(\SL_2(A), \SM_2(A); X_p(A^2), \hat{X}_p(A^2)) & \text{if $p=0,1$}\\
	H_q(\SL_2(A), \SM_2(A); Z_1(A^2), \hat{Z}_1(A^2))  & \text{if $p=2$}\\
	0 & \text{if $p>2$}
\end{cases}
\!\Rightarrow\! H_{p+q}(\SL_2(A), \SM_2(A); \z).
\]
Consider the long exact sequence
\[
\cdots \arr H_q(\SM_2(A), \hat{X}_p(A^2)) \arr H_q(\SL_2(A),X_p(A^2))
\arr \Eb_{p,q}^1 \arr H_{q-1}(\SM_2(A), \hat{X}_p(A^2)) 
\]
\[
\arr H_{q-1}(\SL_2(A),X_p(A^2)) \arr \cdots.
\]
Since
\[
H_q(\SL_2(A), X_0(A^2))\simeq H_q(B(A),\z), \ \ \ H_q(\SL_2(A), X_1(A^2))\simeq H_q(T(A),\z),
\] 
and 
\[
H_q(\SM_2(A), \hat{X}_0(A^2))\simeq H_q(T(A),\z), \ \ \ H_q(\SM_2(A), \hat{X}_1(A^2))\simeq H_q(T(A),\z),
\]
from the above exact sequence, for any $q$, we get 
\[
\Eb_{0,q}^1\simeq \SS_q\simeq H_q(B(A),T(A);\z), \ \ \ \ \Eb_{1,q}^1=0.
\]
Therefore
\[
\Eb_{0,q}^2\simeq \Eb_{0,q}^1, \ \ \ \ \Eb_{1,q}^2=0, \ \ \ \ \Eb_{2,q}^2\simeq\Eb_{2,q}^1.
\]
Now by easy analysis of the spectral sequence we get the exact sequence
\[
\cdots \!\arr\! H_{n+2}(\SL_2(A),\! \SM_2(A); \z) \!\arr\! \Eb_{2,n}^2\!\! \arr\! H_{n+1}(B(A),\!T(A);\z)\! \arr 
\!H_{n+1}(\SL_2(A),\! \SM_2(A); \z)
\]
\[
\arr \Eb_{2,n-1}^2 \arr H_{n}(B(A),T(A);\z)\arr \cdots
\]
where $H_{n}(B(A),T(A);\z) \arr H_{n}(\SL_2(A), \SM_2(A); \z)$ is induced by the natural inclusion
of the pairs $(B(A),T(A)) \harr (\SL_2(A), \SM_2(A))$. 

It is easy to see that $\Eb_{0,0}^2=0=\Eb_{1,0}^2$. Moreover we have the exact sequence
\[
H_0(\SM_2(A), \z) \arr H_0(\SL_2(A), Z_1(A^2)) \arr  \Eb_{2,0}^2 \arr 0.
\]
Note that $H_0(\SM_2(A), \z) \simeq \z$ and $H_0(\SL_2(A),Z_1(A^2))=\GW(A)$. Moreover, the map 
$\z \arr \GW(A)$ is injective and  sends $1$ to $p_{-1}^+=\lan -1\ran +1$. Thus 
\[
\Eb_{2,0}^2 \simeq \GW(A)/\lan \lan -1\ran +1\ran.
\]
Furthermore we have the exact sequence
\[
H_1(\SM_2(A), \z) \arr H_1(\SL_2(A), Z_1(A^2)) \arr  \Eb_{2,1}^2 \arr 0.
\]
From the commutative diagram
\[
\begin{tikzcd}
	H_1(\SM_2(A), \z) \ar[r]\ar[d] & H_1(\SL_2(A),Z_1(A^2)) \ar[r]\ar[d] & \Eb_{2,1}^2 \arr 0\ar[d]\\
	H_1(\SM_2(A), \hat{X}_1(A^2) \ar[r] & H_1(\SL_2(A), X_1(A^2) \ar[r]& 0
\end{tikzcd}
\]
we obtain the exact sequence
\[
\GG_A \overset{\alpha}{\larr} E_{2,1}^2 \arr  \Eb_{2,1}^2 \arr 0.
\]
On the other hand we have the exact sequence
\begin{align*}
	H_3(B(A),T(A);\z) \arr H_3(\SL_2(A), \SM_2(A); \z) \arr \Eb_{2,1}^2 \arr H_2(B(A),T(A);\z) \arr\\
	H_2(\SL_2(A), \SM_2(A); \z) \arr \frac{\GW(A)}{\lan \lan -1\ran +1\ran} \arr A_\aa
	\arr H_1(\SL_2(A), \SM_2(A); \z) \arr 0.
\end{align*}
We call $\GW(A)/\lan \lan -1\ran +1\ran$ the {\it Witt group} of $A$ and denote it with $W(A)$:
\[
W(A):=\GW(A)/\lan \lan -1\ran +1\ran.
\]
This can be seen as a natural generalization of the classical Witt group of a field 
(see \cite[Definition 1.3, Chap II]{lam2005}).

\begin{prp}\label{h2-h3}
	Let $A$ be a $\GE_2$-ring such that $H_i(T(A),\z)\simeq H_i(B(A),\z)$ for $i\leq 3$. Then
	\par {\rm (i)} $H_2(\SL_2(A), \SM_2(A); \z) \simeq  W(A)$.
	\par {\rm (ii)} $H_3(\SL_2(A), \SM_2(A); \z) \simeq \Eb_{2,1}^2$. In particular, we have the exact sequence
	\[
	\GG_A \overset{\alpha}{\larr} E_{2,1}^2 \arr H_3(\SL_2(A), \SM_2(A); \z) \arr 0.
	\]
\end{prp}
\begin{proof}
	It follows from our hypothesis that $H_i(B(A),T(A);\z)=0$ for $0\leq i\leq 3$. Now the claims follows from 
	the above discussions.
\end{proof}

\begin{thm}\label{h3-RP}
	Let $A$ be a universal $\GE_2$-ring such that $H_i(T(A),\z)\simeq H_i(B(A),\z)$ for $i\leq 3$. 
	Then we have an exact sequence
	\[
	I(A)\otimes \mu_2(A) \arr H_3(\SL_2(A), \SM_2(A); \z) \arr \displaystyle\frac{\RP_1(A)}{\lan 
		\psi_1(a^2):a\in \aa\ran} \arr 0.
	\]
	In particular, if $-1\in (\aa)^2$, then $H_3(\SL_2(A),\SM_2(A);\z)\simeq \RP_1(A)$.
\end{thm}
\begin{proof}
	The first claim follows from the above Proposition, Lemma \ref{image} and
	the following diagram with exact row and column:
	\[
	\begin{tikzcd}
		& \II_A\otimes \mu_2(A)\ar[d] & \\
		\GG_A \ar[r, "\alpha"] & E_{2,1}^2 \ar[r] \ar[d]& H_3(\SL_2(A), \SM_2(A); \z)\arr 0\\
		& \RP_1(A)\ar[d] & \\
		& 0& 
	\end{tikzcd}
	\]
	(Note that in above diagram we may replace $\II_A$  with $I(A)$.)
	The second claim follows from the first claim, Lemma \ref{crucial} and the fact that $\psi_1(a^2)=0$ 
	(see \cite[Lemma 3.21(5)]{C-H2022}).
\end{proof}


\begin{thm}\label{witt}
	Let $A$ be ring such that $H_i(T(A),\z)\simeq H_i(B(A),\z)$ for $i\leq 3$. Let $H_1(\SL_2(A),\z)=0$. 
	\par {\rm (i)} If $A$ is a $\GE_2$-ring, then  $H_2(\SL_2(A),T(A);\z)\simeq K_1^\MW(A)$.
	\par {\rm (ii)} If $A$ is a universal $\GE_2$-ring, then $H_3(\SL_2(A),T(A);\z\half)\simeq \RP_1(A)\half$.
\end{thm}
\begin{proof}
	(i) From the inclusions $T(A)\se \SM_2(A)\se \SL_2(A)$, we obtain the long exact sequence
	\[
	\cdots\arr H_n(\SM_2(A), T(A);\z) \arr H_n(\SL_2(A), T(A);\z) \arr H_n(\SL_2(A), \SM_2(A);\z) \arr 
	\]
	\[
	H_{n-1}(\SM_2(A), T(A);\z) \arr \cdots.
	\]
	Since $H_1(\SL_2(A),\z)=0$, we have 
	\[
	H_1(\SL_2(A), \SM_2(A);\z)=0=H_1(\SL_2(A), T(A);\z).
	\]
	It is easy to see that
	\[
	H_1(\SM_2(A), T(A);\z)\simeq \z/2.
	\]
	We already have seen that  $H_2(\SL_2(A), \SM_2(A);\z)\simeq W(A)$ (Proposition \ref{h2-h3}). 
	Form the exact sequences
	\[
	H_2(T(A),\z)\!\arr\! H_2(\SL_2(A),\z) \!\arr\! H_2(\SL_2(A), T(A);\z) \!\arr\! 
	H_1(T(A),\z) \!\arr\! H_1(\SL_2(A),\z)=0
	\]
	and 
	\[
	H_2(T(A),\z)\arr H_2(\SL_2(A),\z) \arr I^2(A) \arr 0
	\]
	we obtain the exact sequence
	\begin{equation}\label{MW1}
		0\arr I^2(A) \arr H_2(\SL_2(A), T(A);\z) \arr K_1^M(A)\arr 0.
	\end{equation}
	Now consider the exact sequence
	\[
	H_2(T(A),\z)\arr H_2(\SM_2(A),\z) \arr H_2(\SM_2(A), T(A);\z) \arr H_1(T(A),\z) 
	\]
	\[
	\arr H_1(\SM_2(A),\z) \arr H_1(\SM_2(A), T(A);\z)  \arr 0.
	\]
	Since $H_2(T(A),\z)\arr H_2(\SM_2(A),\z)$ is surjective (by Lemma \ref{H2SM2}) and 
	$H_1(\SM_2(A),\z)$ sits in the exact sequence
	$1 \arr \GG_A \arr H_1(\SM_2(A),\z)\arr \z/2\arr 0$, we have
	\[
	H_2(\SM_2(A), T(A);\z) \simeq \aa^2\simeq 2K_1^M(A).
	\]
	Thus we have the exact sequence
	\begin{equation}\label{MW2}
		0\arr 2K_1^M(A)\arr H_2(\SL_2(A), T(A);\z) \arr I(A) \arr 0.
	\end{equation}
	
	It is known that the first Milnor-Witt $K$-group of $A$, $K_1^\MW(A)$, satisfies in the 
	exact sequences (\ref{MW1}) and (\ref{MW2}) (\cite[\S 2]{hutchinson-tao2010}). 
	From the exact sequences (\ref{MW1}) and (\ref{MW2}) we obtain the commutative diagram
	\[
	\begin{tikzcd}
		H_2(\SL_2(A),T(A);\z)\ar[r]\ar[d] & K_1^M(A)\ar[d]\\   
		I(A) \ar[r] & I(A)/I^2(A).
	\end{tikzcd}
	\]
	Since $I(A)/I^2(A)\simeq \GG_A\simeq K_1^M(A)/2K_1^M(A)$,  the above diagram is Cartesian. Thus
	\[
	H_2(\SL_2(A),T(A);\z)\simeq K_1^M(A) \times_{I(A)/I^2(A)} I(A).
	\]
	But it is well-known that $K_1^\MW(A)$ is the Cartesian product of the maps $K_1^M(A) \arr I(A)/I^2(A)$
	and $I(A) \arr I(A)/I^2(A)$ (or we can take this as definition). Thus
	\[
	H_2(\SL_2(A),T(A);\z)\simeq K_1^M(A) \times_{I(A)/I^2(A)} I(A) \simeq K_1^\MW(A).
	\]
	(ii) Consider the long exact sequence
	\[
	H_3(\SM_2(A), T(A);\z) \arr H_3(\SL_2(A), T(A);\z) \arr H_3(\SL_2(A), \SM_2(A);\z)
	\]
	\[
	\arr 2K_1^M(A) \arr K_1^\MW(A) \arr W(A).
	\]
	This gives us the exact sequence
	\[
	H_3(\SM_2(A), T(A);\z) \arr H_3(\SL_2(A), T(A);\z) \arr H_3(\SL_2(A), \SM_2(A);\z) \arr 0.
	\]
	Consider the exact sequence 
	\[
	H_3(T(A),\z) \!\arr\! H_3(\SM_2(A),\z) \!\arr\! H_3(\SM_2(A), T(A);\z) \!\arr\! 
	H_2(T(A),\z)\!\arr\! H_2(\SM_2(A),\z)
	\]
	We have seen that the kernel of the right hand side map is isomorphic to $\aa\wedge \mu_2(A)$.
	Moreover using the spectral sequence $\hat{E}_{p,q} \Rightarrow H_{p+q}(\SM_2(A),\z)$ we obtain
	the exact sequence
	\[
	0 \arr (\aa \wedge \aa)/2 \arr H_3(\SM_2(A),\z)/H_3(T(A),\z) \arr \GG_A \arr \aa \wedge \aa.
	\]
	These show that $H_3(\SM_2(A), T(A);\z\half)=0$ Thus
	\[
	\begin{array}{c}
		H_3(\SL_2(A), T(A);\z\half) \simeq H_3(\SL_2(A), \SM_2(A);\z\half)\simeq \RP_1(A)\half.
	\end{array}
	\]
\end{proof}

\begin{rem}
	It is known that $K_1^\MW(A)$ and $\RP_1(A)$ have certain localization property 
	\cite[Theorem~6.3]{{gsz2016}}, \cite[Theorem ~A]{hmm2022}. Wendt in \cite[App. A]{wendt2018} 
	have introduced a higher version of these groups. It would be interesting to
	see what is the connection of these groups to the relative homology groups $H_n(\SL_2(A), \SM_2(A);\z\half)$.
\end{rem}


\end{document}